\documentclass[sn-basic]{sn-jnl}

\usepackage{graphicx}%
\usepackage{multirow}%
\usepackage{amsmath,amssymb,amsfonts}%
\usepackage{amsthm}%
\usepackage{mathrsfs}%
\usepackage[title]{appendix}%
\usepackage{xcolor}%
\usepackage{textcomp}%
\usepackage{manyfoot}%
\usepackage{booktabs}%
\usepackage{algorithm}%
\usepackage{algorithmicx}%
\usepackage{algpseudocode}%
\usepackage{listings}%

\usepackage{caption}
\usepackage{subcaption}
\usepackage{graphicx}%
\usepackage{multirow}%
\usepackage{amsmath,amssymb,amsfonts}%
\usepackage{amsthm}%
\usepackage{mathrsfs}%
\usepackage{latexsym}%
\usepackage{mathtools}%
\usepackage[title]{appendix}%
\usepackage{xcolor}%
\usepackage{color,cancel}
\usepackage{textcomp}%
\usepackage{manyfoot}%
\usepackage[utf8]{inputenc}
\usepackage{listings}%

\usepackage{lipsum}
\usepackage{wrapfig}
\usepackage{enumerate}
\usepackage{soul}
\usepackage{chemarrow}

\usepackage{booktabs}%
\usepackage{multicol}
\usepackage{multirow}
\usepackage[math]{cellspace}
\usepackage{algorithm}%
\usepackage{algorithmicx}%
\usepackage{algpseudocode}%
\usepackage{listings}%
\usepackage{units}
\usepackage{verbatim}
\usepackage{fancyvrb}

\setkeys{Gin}{width=\linewidth,totalheight=\textheight,keepaspectratio}

\usepackage{hyperref}
\hypersetup{
	hidelinks = true
}

\renewcommand{\d}{\mathrm{d}}
\newcommand{\rbb}{\mathbb{R}}

\newcommand{\lambdags}{\lambda_{\mathrm{g} \to \mathrm{s}}}
\newcommand{\lambdasg}{\lambda_{\mathrm{s} \to \mathrm{g}}}

\newcommand{\catform}{\phi}

\newcommand{\mtzero}{\psi}
\newcommand{\Tnd}{\widetilde{T}}

\newcommand{\lcrit}{\ell_\mathrm{crit}}
\newcommand{\Lchar}{L_\mathrm{char}}
\newcommand{\Lcrit}{L_\mathrm{crit}}

\newcommand{\ndlambdamin}{\widetilde{\lambda}_\mathrm{min}}
\newcommand{\tchar}{t_\mathrm{char}}

\newcommand{\micron}{\mu \mathrm{m}}
\newcommand{\domain}{\mathcal{O}}

\newtheorem{theorem}{Theorem}[section]
\newtheorem{proposition}[theorem]{Proposition}
\newtheorem{definition}[theorem]{Definition}

\newtheorem{remark}[theorem]{Remark}

\begin{document}

\title[Minimal Mechanisms of Microtubule Length Regulation in Living Cells]{Minimal Mechanisms of Microtubule Length Regulation in Living Cells}


\author*[1]{\fnm{Anna C} \sur{Nelson}}\email{anelson@math.duke.edu}

\author[2]{\fnm{Melissa M} \sur{Rolls}}

\author[1,3]{\fnm{Maria-Veronica} \sur{Ciocanel}}
\equalcont{These authors contributed equally to this work.}

\author[4]{\fnm{Scott A} \sur{McKinley}}
\equalcont{These authors contributed equally to this work.}

\affil[1]{\orgdiv{Department of Mathematics}, \orgname{Duke University}, \orgaddress{\city{Durham}, \postcode{27710}, \state{NC}, \country{USA}}}

\affil[2]{\orgdiv{Department of Biochemistry and Molecular Biology}, \orgname{Pennsylvania State University}, \orgaddress{ \city{State College}, \postcode{16801}, \state{PA}, \country{USA}}}

\affil[3]{\orgdiv{Department of Biology}, \orgname{Duke University}, \orgaddress{ \city{Durham}, \postcode{27710}, \state{NC}, \country{USA}}}
\affil[4]{\orgdiv{Department of Mathematics}, \orgname{Tulane University}, \orgaddress{ \city{New Orleans}, \postcode{70118}, \state{LA}, \country{USA}}}

\abstract{The microtubule cytoskeleton is responsible for sustained, long-range  intracellular transport of mRNAs, proteins, and organelles in neurons. Neuronal microtubules must be stable enough to ensure reliable transport, but they also undergo dynamic instability, as their plus and minus ends continuously switch between growth and shrinking. This process allows for continuous rebuilding of the cytoskeleton and for flexibility in injury settings. Motivated by \textit{in vivo} experimental data on microtubule behavior in \textit{Drosophila} neurons, we propose a mathematical model of dendritic microtubule dynamics, with a focus on understanding microtubule length, velocity, and state-duration distributions. We find that limitations on microtubule growth phases are needed for realistic dynamics, but the type of limiting mechanism leads to qualitatively different responses to plausible experimental perturbations. We therefore propose and investigate two minimally-complex length-limiting factors: limitation due to resource (tubulin) constraints and limitation due to catastrophe of large-length microtubules. We combine simulations of a detailed stochastic model with steady-state analysis of a mean-field ordinary differential equations model to map out qualitatively distinct parameter regimes. This provides a basis for predicting changes in microtubule dynamics, tubulin allocation, and the turnover rate of tubulin within microtubules in different experimental environments.}

\keywords{microtubule turnover, stochastic modeling, ordinary differential equations, parameterization}

\maketitle
	
\section{Introduction}\label{sec:intro}

Microtubules (MTs) are protein filaments that provide structure to cells and allow them to take on diverse shapes. For instance, neurons are long-range cells where protein transport must occur across long distances, and this is crucially supported by the microtubule cytoskeleton \citep{kapitein2011way,kapitein2015building,kelliher2019microtubule,Rolls2021}. These polymers must therefore be stable enough to support reliable transport of various proteins and organelles. However, microtubules also have dynamic ends, which contribute to the continuous rebuilding of the cytoskeleton. For example, turnover experiments in segments of dendrites of \textit{Drosophila} neurons show that microtubule polymers, comprised of tubulin, turn over every few hours in these cells \citep{Tao2016}. This is due to the switch-like growth behavior that both ends exhibit, where periods of growth transition suddenly into shrinking behavior (catastrophe) and back to growth (rescue). This behavior, called dynamic instability \citep{mitchison1984dynamic}, is especially important when neurons are confronted with stress and injury \citep{Rolls2021}.

Given the complex and intrinsically stochastic dynamics of the microtubule cytoskeleton, many prior mathematical models of MT dynamics have been developed to describe it. These works have used various mathematical frameworks, but usually their primary goal was to reproduce fine details of dynamic instability of MTs observed from \textit{in vitro} experiments. The approaches proposed included two-phase models of microtubule dynamics, where polymers can switch between growing and shrinking phases and the dynamics are primarily described using coupled partial differential equations \citep{Hill1984,Dogterom1993,Govindan2004, akhmanova2022mechanisms,honore2019growth,barlukova2018mathematical}. Follow-up studies proposed stochastic and continuous models of microtubule growth and catastrophe incorporating GTP hydrolysis mechanisms \citep{flyvbjerg1996microtubule,Deymier2005,Margolin2006,Antal2007,hinow2009continuous,Mazilu2010,Buxton2010,Padinhateeri2012,barlukova2018mathematical} and reflecting the observation that catastrophe may be a multi-step process \citep{bowne2013microtubule}. More complex computational models combined known or assumed mechanical, thermodynamic, and kinetic mechanisms for tubulin assembly and bending, GTP hydrolysis, or molecular dynamics of tubulin-tubulin interactions to understand the structure and fluctuation of microtubule tips during dynamic instability \citep{molodtsov2005molecular,Ji2011,margolin2012mechanisms,Zakharov2015}. 

Despite the dynamic instability of individual filaments, microtubule networks are remarkably stable and long-lived, particularly in neurons \citep{Rolls2021}. One population-scale quantity that has attracted attention from mathematical modelers is the distribution of microtubule lengths. While some studies observed and characterized the transition between parameter regimes yielding unbounded growth and finite steady-state MT lengths, fewer models incorporated limited tubulin availability or the impact of tubulin diffusion on MT length control \citep{Dogterom1995,Deymier2005,hinow2009continuous,Margolin2006,Buxton2010}. Even fewer considered mechanisms where MT regulation may arise from the dependence of certain parameters on the length of the MT. For instance, the model proposed in \cite{Mazilu2010} assumed that the attachment and detachment rates of GTP monomer units to the filament depend linearly on the MT length, motivated by the physical constraint that MT lengths cannot be negative, nor exceed a maximum length. In \cite{Rank2018}, the authors included a tubulin resource limitation in their model by allowing the filament polymerization rate to decrease linearly with increasing MT length, while in \cite{bolterauer1999models}, they assumed that the catastrophe probability is an increasing function of length. However, a comprehensive characterization of how multiple mechanisms such as competition for tubulin resources and length-dependent catastrophe of MTs contribute to MT regulation is lacking. 

In addition, a few of these previous works considered realistic living cell behavior, such as regulation of MT length in the mitotic spindle \citep{Dogterom1995,Rank2018} or the role of cell geometry, nucleation center, and cell edges on dynamic instability \citep{Deymier2005,Margolin2006,Buxton2010}. These studies focused on models of MTs growing from a centrosome or nucleation center and thus exhibiting dynamics at the plus end only. Most studies used kinetic parameters for MT dynamics extracted from \textit{in vitro} studies or otherwise simulated the MT behavior for large parameter ranges. However, arrays of non-centrosomal MTs are key for the specialized functions of many cells \citep{keating1999centrosomal,bartolini2006generation,sanchez2017microtubule,akhmanova2022mechanisms}. In many non-centrosomal arrays, minus ends are localized to specific regions by nucleation or capping proteins including the g-tubulin ring complex and CAMSAPs, while the plus end remains free and dynamic. In neurons, MTs form an overlapping array with plus and minus ends staggered throughout axons and dendrites \citep{baas2011hooks,kapitein2015building}. Plus ends grow and shrink throughout axons and dendrites \citep{baas2016stability,Rolls2021}, and at least in \textit{Drosophila} fruit flies and in zebrafish, minus ends also exhibit dynamic instability \citep{Feng2019}.

Addressing the question of how MT length is regulated in a living cell is especially important when studying cytoskeleton dynamics during cell growth or response to injury. While previous complex models of microtubule behavior have provided useful insights for the understanding of dynamic instability \textit{in vitro}, studying whole-cell cytoskeleton dynamics may not require the same level of small-scale biochemical details. Understanding how local fluctuations in MTs impact the overall stability of the MT cytoskeleton in a healthy living cell requires a minimal modeling approach with parameters that can be validated with experiments. An advantage of a minimal modeling strategy is that it would not attempt to incorporate the effects of individual molecules on microtubule growth at tips, which is still actively studied experimentally. Such a modeling framework would also enable further investigation of microtubule dynamics during developmental or injury processes.

Even with such a minimalistic modeling framework, there are numerous challenges. For instance, it is difficult to fully observe data on global microtubule behavior (such as the polarity, density, or length of microtubules) in living cells. Local information on microtubule behavior can be more readily quantified, but may only reflect partial observations. For example, run lengths of microtubule ends may be cut off due to experimental time windows, and parameters extracted from such measurements are likely to represent effective quantities, which already incorporate mechanisms of MT length regulation in an unknown way.  

In summary, the main goal of this work is to investigate two minimal length-limiting factors---limitation due to resource (tubulin) constraints and limitation due to catastrophe of large-length microtubules---and to understand how these different constraints then shape steady-state statistics of microtubule dynamics. Moreover, we study how steady-state length distributions evolve in response to parameter changes. Since the action of these limiting factors is difficult to capture directly, a predictive stochastic model enables the study of a variety of experimentally observable emergent properties. The stochastic modeling framework has the advantage that it can capture fluctuations in behavior within and among different microtubule trajectories, however it requires many parameters in a parameter space that is difficult and computationally costly to navigate. {We therefore propose a coarse-grained deterministic (ordinary differential equations) model that explicitly addresses certain stochastic effects in a novel way. We find that the average MT length is not sufficient to identify all parameters of interest in the stochastic model. However, the deterministic framework helps restrict our study to parameter combinations that generate target MT lengths, thus reducing the feasible parameter space. } 
As a consequence, we can systematically investigate the large parameter space and make faithful predictions in multiple parameter regimes. Informing the stochastic model in this way allows us to distinguish between the contributions of different mechanisms of microtubule length regulation, and to provide predictions for experiments that may soon be possible, such as controlling tubulin concentration \textit{in vivo}.

\subsection{Overview of modeling decisions}\label{sec:intro_overview_modeling}

We propose a stochastic model that tracks the dynamics of several microtubules along one spatial dimension. This one-dimensional space assumption is appropriate for microtubules extending along a segment of an axon or dendrite of a neuron. We model polymerization and depolymerization dynamics at both the plus and minus end of the microtubules. Incorporating the dynamics at both ends is seldom studied in prior models, especially those motivated by mitotic spindle dynamics \citep{Deymier2005,Dogterom1995,Rank2018,bowne2013microtubule}. We are motivated by systems in which microtubules are dynamic at both ends, consistent with recent evidence on dynamic instability at both microtubule ends in neurons of \textit{Drosophila} fruit flies and zebrafish \citep{Feng2019,thyagarajan2022microtubule}. In addition to overall growth and shrinking of microtubules, these systems may also exhibit treadmilling (growth at one end and shortening at the other end) \citep{rodionov1997microtubule,margolis1998microtubule,white2015microtubule}. As motivated above, we are interested in understanding how MT length distributions are regulated by one or two qualitatively different growth control mechanisms: a regime that is limited by competition for tubulin, and a regime that is limited by intrinsic length-driven limitation of catastrophe. 

The assumption that microtubules compete for tubulin resources is widely accepted and has been implemented in previous mathematical models \citep{Dogterom1995,Deymier2005,hinow2009continuous,Margolin2006,Buxton2010}; in our framework, this potential limitation due to tubulin availability impacts the growth velocity of the microtubules at both ends. While we do not include explicit spatial diffusion of tubulin, we consider separate pools of free and unavailable tubulin. Unlike previous models that included a linear dependence of polymerization rates on free tubulin concentration \citep{bolterauer1999models,Rank2018}, we use Michaelis--Menten kinetics to model a smooth dependence of the growth velocity on the free tubulin pool. 

Modeling a length-dependent catastrophe mechanism is motivated by observations from \textit{in vitro} and \textit{in vivo} data \citep{gardner2013microtubule}. Experiments have shown that catastrophe is suppressed in short newly polymerizing microtubules \citep{gardner2011depolymerizing}, and the initial simple hypotheses of stochastic changes in GTP cap length leading to catastrophe could not account for the increase in catastrophe frequency as microtubules grow. Two types of mechanisms have been proposed to account for the increased catastrophe frequency as microtubules grow. The first is the ``antenna" mechanism for kinesin-stimulated catastrophe, which hypothesizes that a depolymerizing kinesin like kinesin-8 drops down on a microtubule and walks to the plus end where it removes subunits to promote catastrophe. Longer microtubules would receive more motors and thus have more subunit removal at the plus end \citep{varga2006yeast,Varga2009}. The second major hypothesized mechanism posits that, as microtubules grow, the plus ends become more ragged, with some protofilaments lagging behind while others keep extending. This accumulation of more variable protofilament lengths as microtubules elongate is supported by observations on plus ends \textit{in vitro} and \textit{in vivo} \citep{coombes2013evolving}, and different modeling approaches support the idea that this could be responsible for length-dependent catastrophe \citep{alexandrova2022theory,coombes2013evolving}. Here, we model length-dependent catastrophe of microtubules by imposing a rate of switching from growth to shrinking that depends linearly on microtubule length.

\begin{figure}
    \centering
    \includegraphics[width=0.8\textwidth]{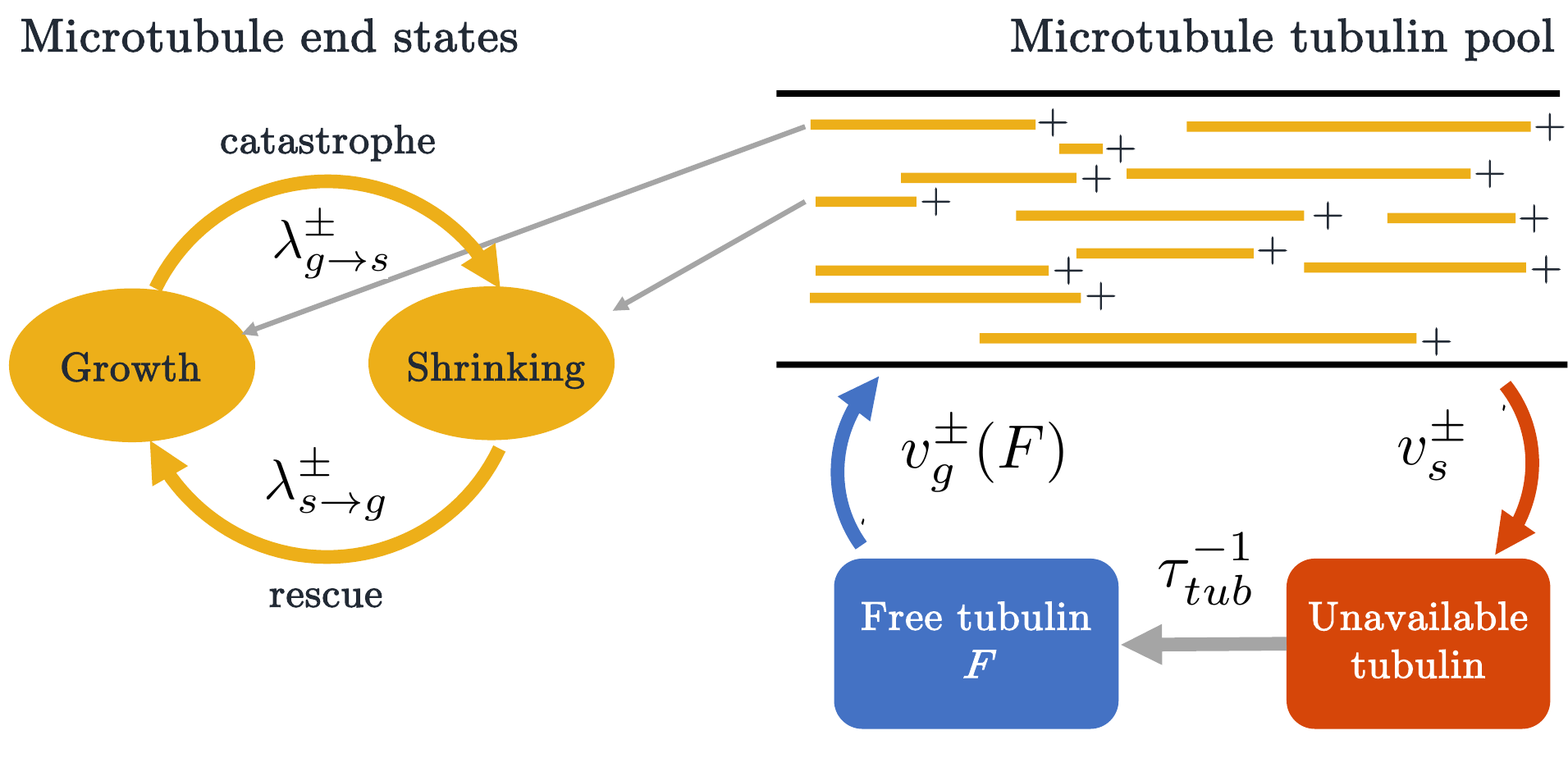}
    \caption{Overview of MT dynamics in a stochastic model simulation with two growth regulation mechanisms: tubulin-dependent growth and length-dependent catastrophe. Each microtubule has a length $L$ and its plus and minus end are in either growth or shrinking states. In growth, the MTs utilize the available tubulin pool; while shrinking, they replenish tubulin into the unavailable pool.  
    }
    \label{fig:ctmc}
\end{figure}

Mathematical models have been developed to probe these and other mechanisms of catastrophe and rescue in microtubule dynamics. Many of these either involved a high level of molecular detail \citep{flyvbjerg1996microtubule,Deymier2005,Margolin2006,Antal2007,hinow2009continuous,Mazilu2010,Buxton2010,Padinhateeri2012,barlukova2018mathematical,bowne2013microtubule,molodtsov2005molecular,Ji2011,margolin2012mechanisms,coombes2013evolving,Zakharov2015}, or incorporated complexity through memory effects, such as age-dependent microtubule catastrophe \citep{barlukova2018mathematical} or location-dependent rescue \citep{fees2019unified}. However, including this level of detail or memory comes with the cost of analytical intractability. Since our goal is to investigate large-scale microtubule behavior in whole cells, and \textit{in vivo} experimental data at this scale is limited, we proceed with  memory-less versions of the switching mechanisms, and adopt a continuous-time Markov chain (CTMC) framework, which is summarized by Figure~\ref{fig:ctmc}. This enables a detailed understanding of qualitative behaviors arising from different parameter regimes and an ability to predict the impact of realistic experimental perturbations.

\subsection{Overview of the paper}
\begin{figure}
    \centering    
    \includegraphics[width=0.9\textwidth]{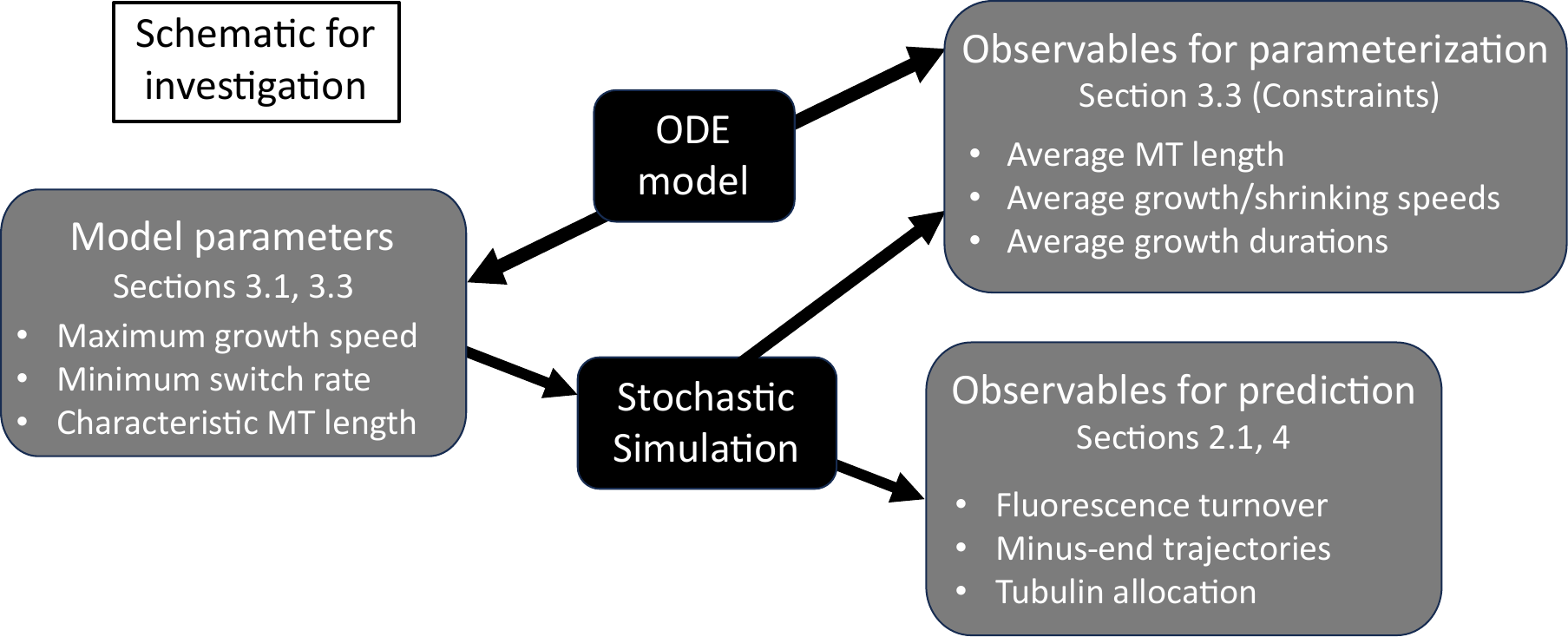}
    \captionsetup{width = 0.9\textwidth}
    \caption{Schematic outlining model parameters of the stochastic simulations, model observables that are used in parameterization of the stochastic model based on the deterministic model, and model outputs that provide predictions for experimental measurements.
    }
    \label{fig:roadmap}
\end{figure}

In Figure~\ref{fig:roadmap}, we show a schematic that outlines our parameterization processes, involving both a stochastic CTMC model and a coarse-grained ordinary differential equations (ODE) model. Section~\ref{sec:stochastic} outlines the CTMC model that describes polymerization and depolymerization of microtubules under various length regulation constraints in a segment of a dendrite. The stochastic model easily relates to experiments, which include observables like minus-end trajectories, average MT length, and growth speed distributions. However, it is difficult to understand the qualitative behavior of the MTs given the large parameter space of this model. Therefore, Section~\ref{sec:ode} introduces the reduced ODE model of MT growth dynamics that we analyze in order to propose a strategy for parameterization of the stochastic model and to identify relevant different parameter regimes. We treat some experimental parameters as more reliably-observed than others, for example, the growth and shrinkage speeds at both MT ends. We use experimental data in \emph{Drosophila} dendrites as our primary guide for parameterization. Then, in order to target certain desired emergent properties like average MT length, we derive a protocol for assigning values to unobservable parameters in Section~\ref{sec:parameterization}. Given this protocol, we then investigate observables for prediction like fluorescence turnover (see Figure~\ref{fig:roadmap}) in Section~\ref{sec:results} based on the parameterized stochastic model in qualitatively different regimes of the MT behavior and under potential experimental manipulations. We discuss our results in Section~\ref{sec:discussion}.

\section{Stochastic model of microtubule dynamics}\label{sec:stochastic}
To investigate the dynamics of growing and shrinking microtubules, we use a stochastic CTMC model. In each simulation, we assume a fixed number $N$ of microtubules. This is similar to previous investigations \citep{Margolin2006,jonasson2020behaviors}, and consistent with the idea that microtubule populations are tightly regulated \textit{in vivo}. We consider a one-dimensional domain, modeling a neuronal dendrite segment that has uniform polarity. This means that each microtubule is summarized by the pair $(x_+(t),x_-(t))$, where $x_+(t)$ and $x_-(t)$ refer to the position of the plus and minus end, respectively. Note that we assume $x_+(t)\ge x_-(t)$ for all $t$. Initially, each MT end is randomly assigned to be in one of two states: growth, where the microtubule end is polymerizing, and shrinking, where the microtubule is depolymerizing. Each microtubule end can switch from one state to another at some rate $\lambda$. We denote the switching rate from growth to shrinking at the plus end of the MT as $\lambda_{g\rightarrow s}^+$, and the rate of switching from shrinking to growth as $\lambda_{s \rightarrow g}^+$. In the shrinking state, the MT plus ends shrink at velocity $v_s^+$, and in the growth state, the plus ends grow with velocity or $v_g^+$.  We use similar notation for describing the dynamics at the minus end. A cartoon depiction of our CTMC stochastic  model framework is shown in Figure \ref{fig:ctmc}.

We consider a small time step $\Delta t = 1$~second for the stochastic model implementation. We then make the assumption that there can be at most one shrinking event in one time step. All possible events that can occur in a time step are shown in Figure~\ref{fig:times} in Appendix~\ref{app:stochmodel}. Depending on the state of the MT end and the time spent in each state in that time step $\Delta t$, we implement shrinking events and then growth events at the MT ends. The MT end positions are then updated according to these events, and we proceed to the next time step. Additional details on the implementation of the stochastic model are provided in Appendix~\ref{app:stochmodel}.

\begin{figure}
     \centering
     \begin{subfigure}[b]{\textwidth}
         \centering
         \includegraphics[width = 0.9\textwidth]{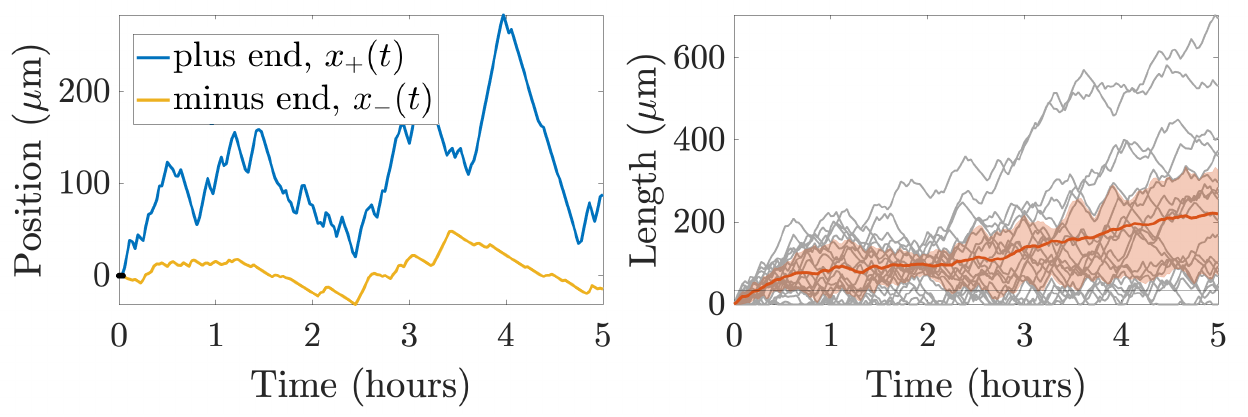}
         \caption{No length-regulating mechanisms}
         \label{fig:sample_nomech}
     \end{subfigure}
     \begin{subfigure}[b]{\textwidth}
         \centering
         \includegraphics[width=0.9\textwidth]{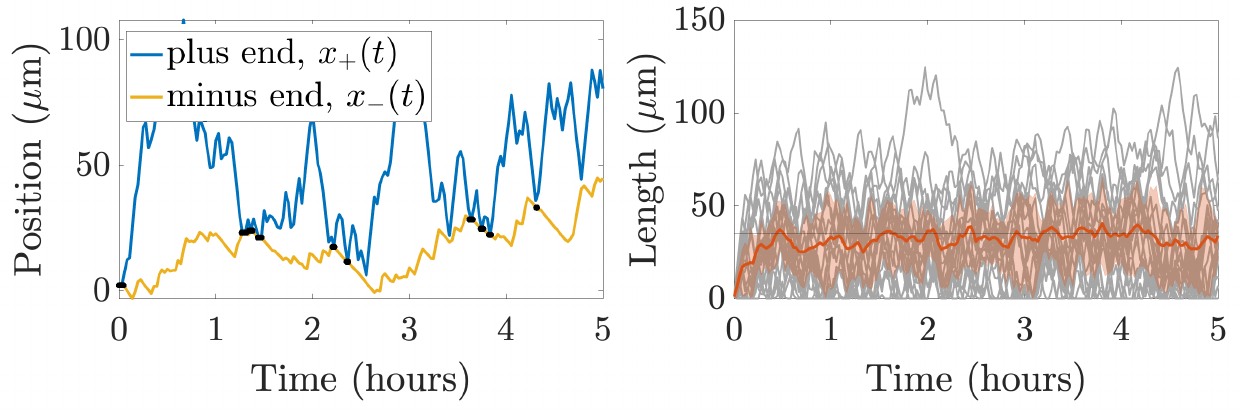}
         \caption{Tubulin-dependent growth and length-dependent catastrophe}
         \label{fig:sample_limtubgamma001}
     \end{subfigure}
    \captionsetup{width = \textwidth}
    \caption{Outputs of a single stochastic simulation for {\bf{(a)}} no length-regulating mechanism and for {\bf{(b)}} length-dependent catastrophe ($\gamma = 0.005$) and tubulin limited growth ($T_{\mathrm{tot}} = 1000\mu$m) for $N=20$ microtubules. The left panel shows positions of the plus end (blue) and minus end (yellow) of a sample MT throughout the entire simulation. Black dots indicate times when the MT reaches length zero. The right panel shows the mean and interquartile range of the MT lengths (orange) and each MT length (grey) throughout the simulation time. The horizonal line denotes the target MT length, $L_* = 35 \mu$m. Note the difference in the $y$-axis range for both panels in (a) and (b). }
        \label{fig:sampletrajectories}
\end{figure}

\textit{In vivo} experimental measurements for plus and minus end growth velocities and run times have been obtained for microtubules in \textit{Drosophila} dendrites in \cite{Feng2019}, while run times and velocities in shrinking states have also been investigated in this system (data unpublished). We can use these experimental quantities to parameterize our CTMC model, allowing us to simulate MT dynamics and observe how filament behavior emerges over time. Figure \ref{fig:sample_nomech} shows trajectory paths for the plus and minus ends of one microtubule (left) and lengths for $20$ microtubules (right) over the course of five hours of simulation time based on the available experimental measurements in \textit{Drosophila} dendrites. The mean microtubule length over time is shown in red, and the red cloud illustrates the interquartile range of lengths over time. The left panel shows that the plus end is moving further to the right through time, while in the right panel, each MT length (grey) is increasing overall throughout time. The length of the dendrites in some \textit{Drosophila} neurons are on the order of tens of microns \citep{Stone2010}, so microtubule lengths on the order of 1000 $\mu$m are not physiologically reasonable for this system. Imposing a domain restriction in this model would lead to microtubules extending throughout the length of the dendrites, which is also unreasonable for this application.

We thus propose two mechanisms to regulate unbounded MT growth: tubulin-dependent growth and length-dependent catastrophe. Figure \ref{fig:ctmc} illustrates how the two mechanisms are incorporated into our CTMC simulation, where growth velocity is a function of free tubulin and the switching from growth to shrinking is dependent on the length of a microtubule ($L$). Specifically, we assume that the switching rate from growth to shrinking satisfies the following conditions.

\begin{definition} \label{defn:catastrophe}
A length-dependent catastrophe rate is defined in terms of four quantities: $L_0 > 0$, the characteristic MT length; $\lambdags^\pm > 0$, the plus-end (respectively, minus-end) growth-to-shrinking phase switch rate when a MT has length $L_0$; $\lambda_{\min} \in (0,\lambdags^\pm)$, the smallest possible value of the catastrophe rate; $\gamma \geq 0$, a parameter that captures the magnitude of length-dependence in the switch rate; and $\catform : \rbb_+ \to \rbb_+$, a dimensionless shape function for the length-dependence.

We say that a function $\catform$ is an admissible length-dependent catastrophe shape function if it is continuous, non-decreasing, and satisfies
\begin{equation} \label{eq:defn-catform}
    \catform(0) = -1, \catform(1) = 0, \text{ and } \catform'(1) = 1.
\end{equation}
When a MT has length $L$, 
\begin{equation}\label{eq:lengthdepcat}
    \text{Length-dependent catastrophe rate: } \max\left(\lambda_{\min}, \lambdags^\pm + \gamma L_0 \phi\left(\frac{L}{L_0}\right)\right).
\end{equation}
\end{definition}
In this way, the catastrophe rate is always strictly positive and equals $\lambdags^\pm$ when $L = L_0$. Moreover, the derivative of the rate equals $\gamma$ at the characteristic length $L_0$. In practice, we will use $\catform(x) = x - 1$ for the shape of the rate function. Note that since $\phi(0) = -1$, when $\gamma < (\lambdags^\pm-\lambda_{\min})/L_0 $ we have that the rate exceeds $\lambda_{\min}$ for all $L$. Otherwise, the minimum rate $\lambda_{\min}$ will be imposed for short MTs.

When including these mechanisms, microtubule growth is facilitated by the availability of tubulin dimers. Tubulin protein can thus either be found in microtubules or be available throughout the neuron for growth. We assume that growth is dependent on the amount of free tubulin $(F)$, and that a shrinking event releases tubulin and thus leads to an increase in $F$. To account for possible slow diffusion of tubulin throughout the length of the dendrite, we assume  that after shrinking, the newly-released tubulin is unavailable for a short period of time, which we denote as unavailable tubulin, $U$. After some time $\tau_{tub}$, the unavailable tubulin $U$ is  made available and becomes free tubulin, $F$. As stated in section \ref{sec:intro_overview_modeling}, growth is modeled using Michaelis--Menten kinetics to yield velocities that saturate in rich tubulin environments and that decrease smoothly to a small growth velocity when free tubulin $F$ is low.

\subsection{Stochastic model observables}\label{sec:model_obs}

Our goal is to use the stochastic model framework outlined in this section to simulate MT dynamics that achieves target measurements such as mean growth speed, mean growth segment duration, and mean MT length (Figure~\ref{fig:roadmap}). Besides these quantities, the model proposed also predicts several emergent measurements that we discuss here. Some of these measurements are quantities that we predict and are not currently available from experiments, while others can be validated with existing experimental work (Figure~\ref{fig:roadmap}). 

For instance, the model quantifies the allocation of tubulin resources into available, unavailable, and microtubule pools. These different tubulin pools are currently very challenging to observe experimentally. In addition, the model can capture a predicted movement direction and speed of the MTs, which emerge as a result of the overall polymerization and depolymerization dynamics at both ends. Figure~\ref{fig:sample_limtubgamma001} shows the predicted trajectory of a single MT in the minus end direction, which suggests some treadmilling of the microtubule.  

The framework proposed can also be used to model and validate the dynamics using photoconversion experiments, which represent an approach for measuring microtubule stability and turnover \textit{in vivo}. 

\begin{figure}
    \centering
\includegraphics[width=0.5\textwidth]{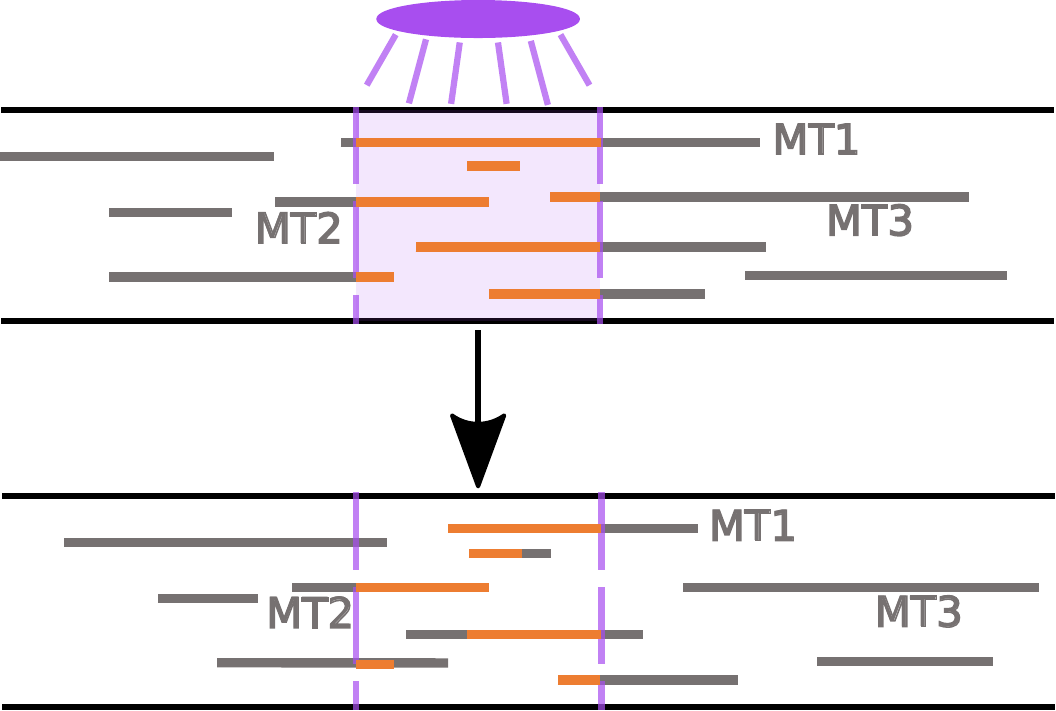}
    \caption{Cartoon of photonvertible tubulin experiment, inspired from \citep{Rolls2021}. {\bf{(Top)}} Tubulin in existing MTs in a window of interest is fluoresced at the photoconversion time. {\bf{(Bottom)}} The intensity of the remaining fluoresced tubulin in the same window is tracked at subsequent times. }
    \label{fig:turnover}
\end{figure}

In axons and dendrites of neurons, photoconversion experiments tag tubulin in a region of interest and observe the fluorescence in that region at subsequent time points \citep{Rolls2021} (see Figure~\ref{fig:turnover}). The fluorescence measured in the window corresponds to tagged tubulin protein that is trapped in a portion of the MT that remains polymerized. Faster decay of the fluorescence in the region of interest indicates more dynamic MTs, which depolymerize on short timescales \citep{Rolls2021}.

We can extract this measurement of MT turnover using our stochastic model as follows. We first select a photoconversion time ($t=120$ minutes), which denotes the time in the simulation when tubulin is tagged in a $10~\mu$m window of interest. We then calculate the microtubule content distribution throughout space at the photoconversion time and choose placement of the window of interest in a way that maximizes the tubulin content (from MTs) in the window. This choice allows us to observe the dynamics of as many microtubule segments as possible in the photoconversion window. At all subsequent times, we track the location of the MT segments that were tagged during the simulated photoconversion and determine the portions of these segments that remain within the window and have not yet depolymerized. As in experiments, we are interested in the relative fluorescence intensity remaining in the window of interest as a function of time \citep{Tao2016,Rolls2021}, and therefore normalize the fluorescence measurement at the photoconversion time to $1$ in all simulation trials. 

For the rest of the paper, we seek to understand what mechanisms affect microtubule behavior, which can be quantified through our stochastic model observables. However, the stochastic simulation relies on many parameters related to growth and shrinking dynamics, and how those parameters depend on each other and influence microtubule growth dynamics is unknown. Additionally, running CTMC simulations with a small timestep is computationally expensive. In the next section, we explore a continuous model of tubulin allocation and microtubule growth to understand how model observables change with different stochastic parameters. Since we ultimately wish to qualitatively match experimental data such as microtubule growth speed and growth duration, we would like to be able to tune our stochastic parameters to match these data.

\section{Reduced deterministic model of microtubule dynamics}\label{sec:ode}

\subsection{Statement of ODE model and justification}

In order to identify and characterize fundamentally distinct parameter regimes, and to establish a protocol for fixing parameters for our investigation, we introduce a deterministic ODE model with a reduced version of the dynamics. Our deterministic approximation is a compartmental model, using the flux of tubulin in and out of the microtubule pool to define the terms of the ODE. To simplify the analysis, we will neglect the unavailable pool of tubulin. In this way, we need only track the amount of tubulin in microtubules at time $t$, $M(t)$, with the amount of free tubulin satisfying $F(t) = T_{\mathrm{tot}} - M(t)$, where $T_{\mathrm{tot}}$ is the (constant) total amount of tubulin in the system. Meanwhile, let $G^+(t)$ and $G^-(t)$ denote the number of MTs that are in growth phase at the plus end and minus end, respectively. We assume a fixed number of MTs, $N$, and let $L(t)=\frac{M(t)}{N}$ denote the average MT length. Similarly, we let $g^+(t)=\frac{G_+(t)}{N}$ and $g^-(t)=\frac{G_-(t)}{N}$ correspond to the fraction of MTs in growth phase at each end.

\begin{table}[h]
\begin{centering}
    \centering
    \captionsetup{width = 0.8\textwidth}
    \begin{footnotesize}
    \begin{tabular}{r|c|cc}
        \hline \hline \textbf{Experimentally-informed quantities} & \textbf{Notation} & \textbf{Plus-end} & \textbf{Minus-end} \\
        \hline Max growth-phase speed & $v_g^{\pm,\mathrm{max}}$ & 9 $\micron/\mathrm{min}$ \citep{Feng2019} & 1.125 $\micron/\mathrm{min}$ \citep{Feng2019} \\
        Average growth-phase speed & $\overline{v}_g^{\pm}$ & 6 $\micron/\mathrm{min}$ \citep{Feng2019} & 0.75 $\micron/\mathrm{min}$ \citep{Feng2019} \\
        Average growth-phase duration & $\overline{\tau}_g^{\pm}$ & 2 min \citep{Feng2019} & 4 min \citep{Feng2019}\\
        Average shrinking-phase speed ($|\text{MT}| > 0$) & $\overline{v}_s^{\pm}$ & 6 $\micron/\mathrm{min}$ $\ddagger$ & 3.5 $\micron/\mathrm{min}$ $\ddagger$ \\
        \hline \hline \textbf{Prescribed model parameters} & \textbf{Notation} & \multicolumn{2}{c}{\textbf{Value(s)}}\\
        \hline Total available tubulin & $T_\mathrm{tot}$ & \multicolumn{2}{c}{$500 \text{ to } 4000 \, \mu\text{m}$} \\
        Number of MT sharing tubulin pool & $N$ & \multicolumn{2}{c}{20} \\
        Average MT length & $\overline{L}$ & \multicolumn{2}{c}{$35 \mu\text{m}$} \\
        Characteristic MT length & $L_0$ & \multicolumn{2}{c}{$35 \mu\text{m}$} \\
        Length-dependence of catastrophe rate & $\gamma$ & \multicolumn{2}{c}{$0 \text{ to } 0.03 \, (\mu\text{m} \cdot \text{min})^{-1}$} \\
        \hline \hline \textbf{Implied model parameters} & \textbf{Notation} & \multicolumn{2}{c}{\textbf{Formula}} \\
        \hline Polymerization/depolymerization speeds &
        $v_g^\pm, \, v_s^\pm$ & \multicolumn{2}{c}{$v_g^{\pm,\max}, \, \overline{v}_s^{\pm}$} \\
        Catastrophe rates at length $L_0$ &
        $\lambda_{g\to s}^\pm$ & \multicolumn{2}{c}{$1/\overline{\tau}_{g}^{\pm}$} \\
        Rescue rates &
        $\lambda_{s\to g}^\pm$ & \multicolumn{2}{c}{see Equation \eqref{eq:defn-lambda-sg}} \\
        Michaelis--Menten constant (scarce tubulin) & $F_{1/2}$ &
        \multicolumn{2}{c}{see Equation \eqref{eq:F-half}} \\
        MT length distribution shape & $\alpha$ & 
        \multicolumn{2}{c}{$1 + \gamma \overline{L} \overline{\tau}_g^+$} \\
        Critical near-zero MT length & $\Lcrit$ & \multicolumn{2}{c}{$v_s^+ \ln 2 / \lambdasg^+$}
        \\ \hline
        \end{tabular}
    \caption{Parameters that appear in the model as well as experimentally observed quantities that inform our choices. Parameters marked with $\ddagger$ are preliminary estimates from the Rolls Lab. A full description of our parameterization strategy is addressed in Section \ref{sec:parameterization}. \label{tab:ode_params}}
    \end{footnotesize}
    \end{centering}
\end{table}

Many of the features of the deterministic model are a direct translation or slight generalization of what we use for stochastic model components. One example is the form of the length-dependent catastrophe, which is here expressed in terms of a dimensionless shape function $\catform(x)$ (see Definition~\ref{eq:defn-catform}). On the other hand, there are aspects of the mean-field model that are not directly represented in the stochastic version. For example, we introduce a function $\mtzero_\alpha(x)$ that phenomenologically captures the tendency of microtubules to hit zero-length during shrinkage runs, using the shape parameter $\alpha$ to capture the effect for different MT length distributions. We express our ODE model using qualitative properties of $\catform$ and $\mtzero_\alpha$, but a specific definition for $\mtzero_\alpha$ can be found in Equation \ref{eq:mtzero-weibull} at the end of our model derivation.

\begin{definition} \label{defn:mtzero}
    Given a shape parameter $\alpha \geq 1$, let the active MT proportion function $\mtzero_\alpha : \rbb_+ \to \rbb_+$ be an increasing and continuously differentiable function satisfying 
    \begin{equation} \label{eq:defn-mtzero}
        \mtzero_\alpha(0) = 0 \text{ and} \lim_{x \to \infty} \mtzero_\alpha(x) = 1. 
    \end{equation}
\end{definition}

\begin{definition} \label{defn:ode}
Let $\catform$ and $\mtzero_\alpha$ satisfy the conditions of Definitions \ref{defn:catastrophe} and \ref{defn:mtzero}, respectively. We define our mean-field approximation of the stochastic model to be the system of ODEs:
    \begin{equation} \label{eq:defn-ode-M}
    \begin{aligned}
        \frac{\d}{\d t} M(t) &= \big(v_g^+ G^+(t) + v_g^- G^-(t)\big) \frac{T_{\mathrm{tot}} - N L(t)}{F_{1/2} + T_{\mathrm{tot}} - N L(t)} \\
        & \qquad \qquad - \Big(v_s^+ (N - G^+) + v_s^- (N - G^-)\Big) \mtzero_\alpha \Big(\frac{L(t)}{\Lcrit}\Big) \\    
        \frac{\d}{\d t} G^{\pm}(t) &= \big(N - G^{\pm}(t)\big) \lambdasg^{\pm} - G^{\pm}(t) \max\left(\lambda_{\mathrm{min}},\lambdags^{\pm} + \gamma L_0 \catform \Big(\frac{L}{L_0}\Big)\right)\,.
    \end{aligned}
    \end{equation}
\end{definition}
The model parameters are described in Table~\ref{tab:ode_params}. By dividing the equations in model \eqref{eq:defn-ode-M} by the constant number of MTs $N$, we obtain the following equations that we focus our subsequent analysis on:
 \begin{equation} \label{eq:defn-ode}
    \begin{aligned}
        \frac{\d}{\d t} L(t) &= \big(v_g^+ g^+(t) + v_g^- g^-(t)\big) \frac{T_{\mathrm{tot}} - N L(t)}{F_{1/2} + T_{\mathrm{tot}} - N L(t)} \\
        & \qquad \qquad - \Big(v_s^+ (1 - g^+) + v_s^- (1 - g^-)\Big) \mtzero_\alpha \Big(\frac{L(t)}{\Lcrit}\Big) \\
        \frac{\d}{\d t} g^{\pm}(t) &= \big(1 - g^{\pm}(t)\big) \lambdasg^{\pm} - g^{\pm}(t) \max\left(\lambda_{\mathrm{min}},\lambdags^{\pm} + \gamma L_0 \catform \Big(\frac{L}{L_0}\Big)\right)\,.
    \end{aligned}
    \end{equation}

In Equations~\eqref{eq:defn-ode}, the fraction of growing microtubules on either the plus  or minus end (denoted as $g^+(t)$ and $g^-(t)$, respectively) is dictated by the flux of microtubules switching into and out of the growing state. The switch from shrinking to growth on either end is given by rate $\lambdasg^{\pm}$, which we assume to be constant. The switch from growth to shrinking is driven by the length-dependent catastrophe mechanism in Definition~\ref{defn:catastrophe}, where the rate of catastrophe is linearly dependent on the MT length $L(t)$, with a minimum rate $\lambda_{\mathrm{min}}$ imposed for short MTs. 

The flux of tubulin into the microtubule pool in Equations~\eqref{eq:defn-ode} is governed by the growth rate of microtubules at the plus and minus ends. We denote the maximum growth speeds by $v_g^+$ and $v_g^-$. However, the actual growth rate can be tempered by the scarcity of available tubulin, which we model with the Michaelis--Menten term in Equations~\eqref{eq:defn-ode}, governed by the half-rate parameter $F_{1/2}$. In Equations~\eqref{eq:defn-ode-M}, before the normalization by the number of MTs $N$, this corresponds to:
\begin{equation} \label{eq:tubulin-flux}
    \text{Tubulin flux into MT pool:} \quad \Big(v_g^+ G^+(t) + v_g^- G^-(t)\Big) \frac{F(t)}{F_{1/2} + F(t)}\,.
\end{equation}

The flux out of the MT pool is difficult to directly model. We assume that in reality the shrinkage rate has a constant mean $v_g^-$ and that this is maintained until the MT hits zero length. In our model, the ``shrinking state'' continues until a switch to growth state occurs (all at Markovian rate $\lambdasg^{\pm}$). The deterministic model is an average over all microtubules, which is to say a weighted average of MTs with non-zero length (with shrinking velocity $v_s^{\pm}$) and of MTs with zero length (hence with shrinking velocity zero). We thus need a term that models the reduction in overall tubulin loss when some MTs in shrinking state have length zero.

Indeed, preliminary simulation work showed that when there is no MT length dependence of the catastrophe rate, the MT length distribution looks exponentially distributed (Figure~\ref{fig:dashboard_baseline_gammas}), and certainly strictly decreasing with length. On the other hand, when MT-length dependence is strong, the distribution is no longer monotone: with a mode that is closer to the median, but is nevertheless quite wide. Throughout our simulations, when the average MT is within a few multiples of the average run length of a shrinking state, then we expect simulations to contain  a non-trivial number of shrinking-state MTs with zero length. 

To capture this translation from an individual-based stochastic model to a mean-field ODE dynamical system, we introduce a family of dimensionless functions $\mtzero_\alpha: \rbb_+ \to [0,1]$, indexed by a \emph{shape parameter} $\alpha > 0$, that we use to approximate the proportion of MTs that have non-zero length when the mean of the length distribution is $L$.  Let $\Lcrit$ be a characteristic length, chosen below, at which a microtubule is ``at risk'' of becoming zero length. Then if $Y$ is a random length drawn from the MT-length distribution with shape $\alpha$ and mean $L$, we write
\begin{equation} \label{eq:psi-alpha}
    \psi_\alpha(L/\Lcrit) \coloneqq P_\alpha(Y > \Lcrit).
\end{equation}

On the right-hand side, we see that this is the probabilistic complement of the cumulative distribution function (cdf) of $Y$. When the distribution of $Y$ is a member of a \emph{shape-and-scale} family, we can express its cdf in terms of the cdf of a unit-scale member of the family that has the same shape parameter. Suppose that $X$ has this unit-scale distribution, and write $F_\alpha(x) = P_\alpha(X \leq x)$. Commonly the unit-scale distribution will have a mean that is different than one, so we define $\mu_\alpha = E_{\alpha}(X).$ Then we can define $Y = \frac{L X}{\mu_\alpha}$ and observe that
\begin{equation}
    P_\alpha(Y > \Lcrit) = P_\alpha\Big(\frac{L X}{\mu_\alpha} > \Lcrit\Big) = 1 - F_\alpha\Big(\frac{\Lcrit \mu_\alpha}{L}\Big)
\end{equation}
In light of \eqref{eq:psi-alpha} this means that we can express $\mtzero_\alpha$ in terms of the unit-scale cdf model:
\begin{equation}
    \mtzero_\alpha(x) := 1 - F_\alpha\Big(\frac{\mu_\alpha}{x}\Big)\,,
\end{equation}
where $x=\frac{L}{\Lcrit}$.
Note that, by the properties of cumulative distribution functions for strictly positive random variables, the conditions of Definition \ref{defn:ode} for $\psi_\alpha$ will commonly be met.

It remains to choose $\Lcrit$ and a family of probability distributions that emulate the qualitative behavior of the MT-length distributions. Since $\Lcrit$ represents a length for which the MT is ``at risk'' of being zero length, we choose to let $\Lcrit$ be the length for which there is a 50\% chance that a MT will hit zero length before the end of a shrinking run. Since the plus-end runs are exponentially distributed with mean $v_s^+/\lambdasg^+$, an appropriate choice for the critical length is the median of that distribution, so
\begin{equation} \label{eq:defn-Lcrit}
    \Lcrit := \frac{v_s^+ \ln 2}{\lambdasg^+}.
\end{equation}
As discussed above, preliminary numerical investigation revealed that when there is no length-dependence in the catastrophe rate, the distribution of MT lengths is qualitatively exponential. As length-dependence becomes more prominent, the distribution is no longer monotone decreasing, but retains very wide support. A natural generalization of the Exponential family is the Weibull family of distributions. This family was originally introduced to model hazard rates (in our context, state-switch rates) that are increasing with time (in our context, with MT length). Moreover, the complement to the cdf has a 
very simple expression: $1 - F_\alpha(x) = \exp(-x^\alpha)$. Since the mean of a $\text{Weibull}(\alpha,1)$ random variable is $\mu_\alpha = \Gamma(1 + \frac{1}{\alpha})$, we have the following final definition for $\mtzero_\alpha$:
\begin{equation} \label{eq:mtzero-weibull}
    \mtzero_\alpha(x) \coloneqq \exp\left(-\left(\frac{\Gamma(1 + \frac{1}{\alpha})}{x}\right)^{\alpha}\right)\,.
\end{equation}

\subsection{Non-dimensionalization and basic mathematical properties}

There are several natural candidates for time- and length-scales to use for the purpose of nondimensionalization. A few important length scales are the average lengths of growth and shrinkage runs (the ratio $v_g^+ / \lambdags^+$ is the average run length of a plus-end growth phase, for example); the total amount of tubulin that is available per microtubule ($T_{\mathrm{tot}}/N$); and the characteristic length $L_0$. We choose the latter, in part because then all quantities can be understood as a fraction of the characteristic MT length scale; then the nondimensionalized dynamics for MT length will be restricted to the interval $[0,T_{\mathrm{tot}}/(N L_0)]$. For the time scale, we choose the average run duration of a plus-end shrinking run ($1/\lambdasg^+$). As we have seen, this quantity plays an important role in establishing the critical length scale of the mean-field approximation for the number of MTs that have non-zero length.
    
\begin{table}[h]
    \centering
    \begin{small}
    \begin{tabular}{c|c}
        Time & Length \\
        \hline 
        $\tchar = 1/\lambdasg^+$ & 
        $\Lchar = L_0$ 
    \end{tabular}

    \vspace{0.5cm}

    \begin{tabular}{c|c|c|c|c|c|c}
        $v^{\pm}$ & $v_g$ & $v_s$ &
        $\lambda^{\pm}$ & $\lambdags$ & $\lambdasg$ &$\ndlambdamin$\\
        \hline \hline
        $\displaystyle \frac{v_g^\pm}{v_s^\pm}$ 
        & $\displaystyle \frac{v_g^+}{v_g^-}$ 
        & $\displaystyle \frac{v_s^+}{v_s^-}$ 
        & $\displaystyle \frac{\lambdags^\pm}{\lambdasg^\pm}$ 
        & $\displaystyle \frac{\lambdags^+}{\lambdasg^-}$ 
        & $\displaystyle \frac{\lambdasg^+}{\lambdasg^-}$ 
        & $\displaystyle \frac{\lambda_{\mathrm{min}}}{\lambdasg^+}$ 
    \end{tabular}

        \vspace*{0.5cm}

    \begin{tabular}{c|c|c|c|c|c|c|c}
        $\Tnd$ & $\eta^\pm$ & $\alpha$ & $\delta_g^\pm$ & $\delta_s^\pm$ & $f_{1/2}$ & $\ell_{\mathrm{crit}}$ & $\ell_0$ \\
        \hline \hline
        $\displaystyle \frac{T_{\mathrm{tot}}}{\Lchar}$ &
        $\displaystyle \frac{\gamma L_0}{\lambdags^\pm}$ 
        & $\displaystyle 1 + \eta^+$
        & $\displaystyle \frac{v_g^\pm / \lambdags^\pm}{\Lchar}$ 
        & $\displaystyle \frac{v_s^\pm / \lambdasg^\pm}{\Lchar}$
        & $\displaystyle \frac{F_{1/2}}{\Lchar}$
        & $\displaystyle \frac{\Lcrit}{\Lchar}$
        & $\displaystyle \frac{L_0}{\Lchar}$
    \end{tabular}

    \caption{Dimensionless groups that appear in the non-dimensional ODE system \eqref{eq:defn-ode-nondim}.}
    \label{tab:nondim}
    \end{small}
\end{table}

Given the choices of length and time scale presented in Table \ref{tab:nondim}, with $\tilde{t}=\frac{t}{\tchar}$ and $\ell=\frac{L}{\Lchar}$, we have the following nondimensional system:
\begin{equation} \label{eq:defn-ode-nondim}
\begin{aligned}
    \ell' &= \delta_s^+ \left(v^+ \Big(g^+ + \frac{g^-}{v_g}\Big) \frac{\Tnd - N \ell}{f_{1/2} + \Tnd - N \ell} - \Big((1 - g^+) + \frac{1}{v_s} (1-g^-)\Big) \mtzero_\alpha\Big(\frac{\ell}{\lcrit}\Big)\right) \\
    (g^{+})' &= (1-g^+) - g^+ \max\left(\ndlambdamin,\Big(1 + \eta^+ \catform \Big(\frac{\ell}{\ell_0}\Big)\Big)\lambda^+\right) \\
    (g^{-})' &= \frac{1}{\lambdasg} \left((1-g^-) - g^- \max\Big(\ndlambdamin,\Big(1 + \eta^- \catform \Big(\frac{\ell}{\ell_0}\Big)\Big)\lambda^- \right)
\end{aligned}
\end{equation}
where $'$ denotes the derivative with respect to the nondimensional time $\tilde t$.

\begin{proposition} \label{prop:ode-well-posed}
Suppose that $\eta^\pm \geq 0$, $\Tnd > N \ell_0$, and all other nondimensional constants in Table \ref{tab:nondim} are positive. Let $\domain = (0,\tilde{T}/N) \times(0,1)\times(0,1)$. Then for any initial condition $(\ell,g^+,g^-)$ in $\domain$ there exists a unique global solution to the system \eqref{eq:defn-ode-nondim}. Moreover, there exists a unique steady-state solution $(\ell_*, g^+_*, g^-_*)$ in the region $\domain$.
\end{proposition}
\begin{remark}
    We do not have a proof that the steady state is asymptotically stable, or that chaos can be excluded from this three-dimensional system. However, numerical experimentation strongly indicates there is convergence to the unique steady-state.
\end{remark}
\begin{proof}
First we observe that the unit cube (the domain closure, $\overline{\domain}$) is a forward-invariant domain for the dynamics. That is to say, the flow from all points on the boundary are to the interior. For example, if either $g^+ = 0$ or $g^- = 0$, then their derivatives are $1$ and $\frac{1}{\lambdasg}$, respectively, and the flow is into the interior of the domain. On the other hand, suppose that either $g^+ = 1$ or $g^- = 1$, the positivity of $\ndlambdamin$ assures that their derivatives are negative. Meanwhile, $0 < g^+, g^- < 1$ and $\ell = 0$ results in $\ell'>0$. Similarly, $\ell = \Tnd/N$ results in $\ell'<0$. To summarize, if the system takes its value in the interior or edge of any faces of the parallelepiped domain, the flow is to the interior of the domain. Finally, we note that none of the corners are fixed points, since the local dynamics pushes the system to a domain face, from which it enters the interior of the domain. For example, at $(0,0,0)$, we have that $\ell' = 0$, but $(g^+)' > 0$ and $(g^-)' > 0$. The domain is therefore a forward invariant set. Because all functions in the system are continuous in the domain with bounded derivatives, they are Lipschitz continuous. Global existence and uniqueness of solutions follows. 

To establish the existence and uniqueness of a steady state, note first that by setting derivatives equal to zero, we can immediately solve for $g^{\pm}_*$, if $\ell_*$ is given. In fact, for any value $\ell$, the growth versus shrinking phase fraction is in balance at the values $g^{\pm}_*(\ell)$ where
\begin{equation} \label{eq:g_star}
    g^\pm_*(\ell) = \frac{1}{1 +  \max\Big(\ndlambdamin,\big(1 + \eta^\pm \catform(\ell/\ell_0)\big)\lambda^\pm\Big)}.
\end{equation}
These values lie strictly between 0 and 1. Moreover, since $\phi(\ell)$ is non-decreasing in $\ell$, both $g^+_*$ and $g^-_*$ are decreasing in $\ell$.

Turning our attention to the equation for $\ell_*$ itself, the steady-state value must satisfy
\begin{equation} \label{eq:ell-star-balance}
   h(\ell_*) =  h_1(\ell_*) - h_0(\ell_*) = 0
\end{equation}
where
\begin{equation} \label{eq:ell-star-parts}
\begin{aligned}
    h_1(\ell) &\coloneqq v^+ \Big(g_*^+(\ell) + \frac{g_*^-(\ell)}{v_g}\Big) \frac{\Tnd - N \ell}{f_{1/2} + \Tnd - N \ell} \,, \\
    h_0(\ell) &\coloneqq \Big((1 - g_*^+(\ell)) + \frac{1}{v_s} (1-g_*^-(\ell))\Big) \mtzero_\alpha\Big(\frac{\ell}{\lcrit}\Big)\,.
\end{aligned}
\end{equation}
Since $\psi_\alpha(0) = 0$, $h_0(0) = 0$. As $\ell$ increases, we note that $h_0(\ell)$ is a product of positive increasing functions, and is hence an increasing function itself. Meanwhile, $h_1(0) > 0$ and $h_1(\Tnd/N) = 0$ while decreasing in $\ell$ in between. It follows that $h(l)$ is decreasing with $h(0)>0$ and $h(\tilde{T}/N)<0$. By the Intermediate Value Theorem, there exists a unique $\ell_* \in (0,\tilde{T}/N)$ such that $h(\ell_*)=0$, so that $\eqref{eq:ell-star-balance}$ holds. It follows that $(\ell_*, g^+_*(\ell_*), g^-_*(\ell_*))$ is a unique fixed point of the system \eqref{eq:defn-ode-nondim}.
\end{proof}

\begin{figure}[h]
    \centering
    \captionsetup{width = 0.9\textwidth}
    \includegraphics[width = 0.9\textwidth]{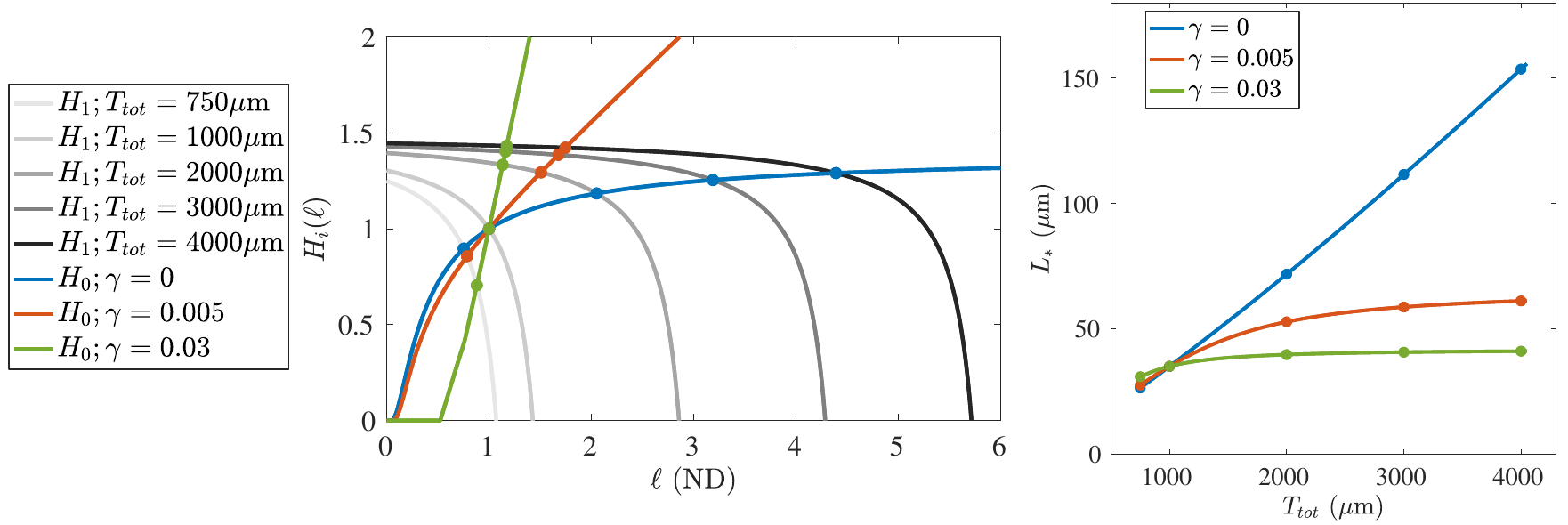}
    \caption{\textbf{(Left)} Various evaluations of the functions $H_0(\ell)$ and $H_1(\ell)$, see Equation \eqref{eq:defn-H}. We use the baseline parameter set in both cases, with the modifications that $H_0$ is evaluated over a range of $\gamma$ (blue, red, and green), while the function $H_1$ is evaluated over a range of total tubulin (light gray to black as $T_{\mathrm{tot}}$ increases.) The steady state values (filled circles) occur at intersections of $H_0$ and $H_1$. \textbf{(Right)} The steady state $L_*$ is displayed as a function of $T_{\mathrm{tot}}$ for the fixed values of $\gamma$. The displayed steady state values (filled circles) correspond to those displayed on the left.}
    \label{fig:H-L-star}
\end{figure}

We display a visualization of the fixed point argument in Figure \ref{fig:H-L-star} through modified versions of the functions $h_0$ and $h_1$. Since the size of the growth phase population is never zero after $t > 0$, we can divide Equation \ref{eq:ell-star-balance} through by $g_*^+(\ell)+ \frac{g_*^-(\ell)}{v_g}$ and define the functions 
\begin{equation} \label{eq:defn-H}
\begin{aligned}
    H_1(\ell) &\coloneqq v^+ \frac{\Tnd - N \ell}{f_{1/2} + \Tnd - N \ell} \,, \\
    H_0(\ell) &\coloneqq \frac{(1 - g_*^+(\ell)) + \frac{1}{v_s} (1-g_*^-(\ell))}{g_*^+(\ell) + \frac{g_*^-(\ell)}{v_g}} \mtzero_\alpha\Big(\frac{\ell}{\lcrit}\Big)\,.
\end{aligned}
\end{equation}
The benefit of this presentation is that the functions $H_1$ and $H_0$ split the dependence on the two major parameters we will perturb in our analysis. The function $H_1$ depends on the total tubulin in the system $T_{\mathrm{tot}}$ through $\Tnd$, but not on the magnitude $\gamma$ of length-dependence in the catastrophe rate. On the other hand, since $\Lchar$ is chosen to be independent of the total tubulin, then $H_0$ depends on $\gamma$ (in non-dimensional form, on $\eta$) through the $g_*^\pm$ functions, but is not affected by changes in $T_{\mathrm{tot}}$. Therefore, as we perturb these two parameters, they lead to a shift in one function or the other, which we can visualize as in the left panel of Figure \ref{fig:H-L-star}. For each combination of $T_{\mathrm{tot}}$ and $\gamma$, there is a unique intersection of $H_0$ and $H_1$, which is known to exist due to the arguments concerning $h_0$ and $h_1$ in the proof of Proposition \ref{prop:ode-well-posed}. 

\subsection{Strategy for parameterization}\label{sec:parameterization}

One of our fundamental goals is to be able to observe a system, find a best-fit set of parameters, and then perturb selected model parameters to predict how the system will respond if such a biological experiment was carried out. An example of such an exercise can be seen in the right panel of Figure \ref{fig:H-L-star}. First, using the procedure described below, we  selected parameters in such a way that the ``experimentally-informed quantities'' listed in Table \ref{tab:ode_params} were respected while setting $L_* = 35\mu m$ when $N = 20$ and $T_\mathrm{tot} = 1000 \mu $m for each choice of $\gamma \in \{0, 0.005, 0.03\}$. The success of the presented procedure can be seen in the left panel of Figure \ref{fig:H-L-star}, because the functions $H_0$ and $H_1$ all intersect at the same location for this fixed value of $T_\mathrm{tot}$. Fixing all parameters except for $T_\mathrm{tot}$, we then solved for $L_*$ as we varied $T_\mathrm{tot}$ from 750 to 4000~$\mu\text{m}$. We display the relationship between $L_*$ and $T_{\mathrm{tot}}$ for three choices of $\gamma$. When $\gamma = 0$, there is no length-dependence in the catastrophe rate, and therefore no penalty in MTs growing larger.  The steady-state average MT length value increases linearly without bound. On the other hand, strong length-dependence in the catastrophe rate results in very little increase in $L_*$ when more tubulin is available in the system. Because the values $\gamma = 0, 0.005, 0.03$ yield qualitatively different responses in $L_*$ to $T_\mathrm{tot}$, we adopt these as our standard choices when conducting numerical experiments for the stochastic model in Section \ref{sec:results}.

In this section, we outline a strategy for selecting model parameters for the stochastic system relying on quantities that are most directly observable. For parameters that are not known, we can derive quantities that will result in fixed points of the ODE that generate the desired values (see Figure~\ref{fig:roadmap}). While not a perfect match to the full stochastic model, we have found very good agreement when the total tubulin in the system $T_\mathrm{tot}$ exceeds the amount required for the MTs to reach their targeted average length.

Referencing the ``experimentally-informed quantities'' and the ``prescribed model parameters'' in Table \ref{tab:ode_params}, we select our parameters so that the following constraints are met:
\begin{enumerate}[(C1)]
    \item the resulting average microtubule length in steady state is $L_* = \overline{L}$; 
    \item the average speed of growth phase polymerization is a fraction of the maximum polymerization speeds in a ratio consistent with the observable quantities $\overline{v}_g^\pm$ and $v_g^{\pm,\max}$; and, 
    \item the growth to shrinking phase switch rate at steady-state matches the inverse of the observed growth phase durations $\overline{\tau}_g^{\pm}$.
\end{enumerate}

Constraints (C2) and (C3) address speeds and durations of MTs in growth phase, however we do not articulate similar constraints for shrinking states. For shrinking state speed, this is because we observed, in preliminary data, that the variance of shrinking phase speeds was significantly less than that of growth phase speeds. We hypothesize that this is because growth phase polymerization is sensitive to tubulin availability, while shrinking phase depolymerization is not subject to an external constraint. We can therefore simply set the model parameters $v_s^\pm$ equal to the experimentally observable parameters $\overline{v}_s^{\pm}$. 

On the other hand, there is not a directly measureable analogue of the switch rate from shrinking to growth. This is because we experimentally measure the time spent in the shrinking state, but there are two reasons shrinking phase runs might stop: either the MT switched back to a growth phase, or the MT reached zero length. This latter possibility is incompatible with the idea that a Markovian switch rate should be the inverse of the average time spent in that state. Effectively, in our model, there is time spent in the ``shrinking state'' while stuck at length zero, causing a mismatch between the model and experimental definitions of the term ``time spent in the shrinking state.'' As a result, we reserve the quantities $\lambdasg^\pm$ as degrees of freedom for the system, that we may use to satisfy other constraints. As we will see, special choices will be made to satisfy (C1).

By construction of the length-dependent catastrophe rate (Definition \ref{defn:catastrophe}), constraint (C3) is immediately satisfied by selecting $L_0 = L_*$ and setting $\lambdags^{\pm} = 1/\overline{\tau}_g^\pm$. Meanwhile, the selection $L_* = \overline{L}$ addresses constraint (C2) if we make a complementary choice for the Michaelis--Menten constant $F_{1/2}$. Indeed, the steady-state flux of tubulin into the microtubule pool can be expressed in two ways: first, $G^\pm$ times the average growth rate of individual MTs that are in each growth phase, i.e. in terms of the experimentally-observable parameters $\overline{v}_g^\pm$; and second, recalling Equation \eqref{eq:tubulin-flux}, in terms of the Michaelis--Menten uptake of tubulin, i.e. in terms of the steady-state value for available tubulin, $F_* = T_\mathrm{tot} - N L_*$ and the maximum MT growth rates. Together, we have the equation 
\begin{equation} \label{eq:tubulin-flux-ss}
    \overline{v}_g^+ G_*^+ + \overline{v}_g^- G_*^- = \left(v_g^{+,\max} G_*^+ + v_g^{-,\max} G_*^-\right) \frac{T_\mathrm{tot} - N L_*}{F_{1/2} + T_\mathrm{tot} - N L_*}\,.
\end{equation}
Now, since the ratios $\overline{v}_g^+/v_g^{+,\max}$ and $\overline{v}_g^-/v_g^{-,\max}$ are hypothesized to be the same (which is equivalent to assuming that the plus- and minus-end Michaelis--Menten constant is the same), we choose to separate the plus- and minus- end constraints and observe that
\begin{displaymath}
    \frac{\overline{v}_g^{\pm}}{v_g^{\pm,\max}} = \frac{T_\mathrm{tot} - N L_*}{F_{1/2} + T_\mathrm{tot} - N L_*}\,,
\end{displaymath}
from which it follows (after imposing the target $L_* = \overline{L}$) that
\begin{equation} \label{eq:F-half}
    F_{1/2} = \left(\frac{v_g^{\pm,\max}}{\overline{v}_g^{\pm}} - 1 \right) \big(T_{\mathrm{tot}} - N \overline{L}\big).
\end{equation}
With this selection, (C2) will be satisfied.

While constraint (C1) may be the most obvious requirement to impose, it is the most difficult to address. By Proposition \ref{prop:ode-well-posed}, we know that there is a unique steady-state solution, so if the number of degrees of freedom match the number of conditions implied at steady-state, then a solution should exist in principle. By virtue of the selection that $L_0 = L_*$, the dimensional version of \eqref{eq:g_star} implies that steady-state values for the proportion of MTs in growth phase at the plus- and minus-ends are
\begin{displaymath}
    G_*^\pm = \frac{\lambdasg^\pm}{\lambdags^\pm + \lambdasg^\pm}.
\end{displaymath}
Substituting these values in the left-hand side of \eqref{eq:tubulin-flux-ss} and setting them equal to the tubulin flux out of the MT pool, $L_*$ must satisfy
\begin{displaymath}
    \frac{\overline{v}_g^+ \lambdasg^+}{\lambdags^+ + \lambdasg^+} + \frac{\overline{v}_g^- \lambdasg^-}{\lambdags^- + \lambdasg^-} = \Big(\frac{\overline{v}_s^+ \lambdags^+}{\lambdags^+ + \lambdasg^+} + \frac{\overline{v}_s^- \lambdags^-}{\lambdags^- + \lambdasg^-}\Big) \mtzero_\alpha \Big(\frac{L_*}{\Lcrit}\Big) \,.
\end{displaymath}
If we separate the plus- and minus-end constraints, recall the definitions of $\mtzero_\alpha$ in Equation \eqref{eq:mtzero-weibull} and $\Lcrit$ in Equation \eqref{eq:defn-Lcrit}, and simplify, we have the system of equations
\begin{equation} \label{eq:lambdasg-almost-there}
    \frac{\overline{v}_g^\pm \lambdasg^\pm}{\overline{v}_s^\pm \lambdags^\pm} = \exp\left(-\left(\frac{\Gamma(1 + \frac{1}{\alpha})  \overline{v}_s^+ \ln 2}{L_* \lambdasg^+ }\right)^\alpha\right) \,.
\end{equation}
This equation can be written more efficiently in terms of the dimensionless quantities $\delta_s^\pm := v_s^\pm / \lambdasg^\pm$. Imposing the choice $\lambdags^\pm = 1/\overline{\tau}_g^\pm$, we can rewrite this as 
\begin{equation} \label{eq:delta-splus-lambert-form}
    \frac{\overline{v}_g^\pm \overline{\tau}_g^\pm}{\Lchar} = \delta_s^\pm \exp\left(-\left(\frac{\Gamma(1 + \frac{1}{\alpha}) \ln 2}{\ell_*} \, \delta_s^+\right)^\alpha \right).
\end{equation}
Looking first at the plus-end version of this equation, $\delta_s^+$ can be expressed in terms of a generalization of the Lambert-$W$ function:
\begin{equation} \label{eq:delta-splus}
    \delta_s^+ = \frac{\ell_*}{\Gamma(1 + \frac{1}{\alpha}) \ln 2} W_\alpha\left(\frac{\overline{v}_g^\pm \overline{\tau}_g^\pm}{L_*}\Gamma(1 + \frac{1}{\alpha}) \ln 2 \right),
\end{equation}
where $W_\alpha(z)$ is the solution of the equation $\displaystyle w e^{-w^\alpha} = z$.

Using the right-hand side of \eqref{eq:delta-splus-lambert-form} to find an expression for the exponential term, we can substitute into \eqref{eq:lambdasg-almost-there} to find a relationship between $\delta_s^-$ and $\delta_s^+$. Namely,
\begin{equation}
    \delta_s^- = \delta_s^+ \frac{\overline{v}_g^- \overline{\tau}_g^-}{\overline{v}_g^+ \overline{\tau}_g^+}. 
\end{equation}
Since the velocity components of $\delta_s^\pm$ are already prescribed, we can satisfy constraint (C1) by choosing $L_* = \overline{L}$ and defining
\begin{equation} \label{eq:defn-lambda-sg}
    \lambdasg^+ \coloneqq \frac{\overline{v}_s^+}{\overline{L}} \frac{\Gamma(1 + \frac{1}{\alpha}) \ln 2}{W_\alpha\left(\frac{\overline{v}_g^\pm \overline{\tau}_g^\pm}{\overline{L}}\Gamma(1 + \frac{1}{\alpha}) \ln 2 \right)}
    \,\, \text{ and }
    \lambdasg^- \coloneqq \lambdasg^+ \frac{\overline{v}_s^-}{\overline{v}_s^+}
    \frac{\overline{v}_g^+ \overline{\tau}_g^+}{\overline{v}_g^- \overline{\tau}_g^-}.
\end{equation}
We have therefore determined all implied model parameters in Table~\ref{tab:ode_params} for our application to MT turnover in \textit{Drosophila} dendrites. This completes our strategy for selecting model parameters for the stochastic model (see Figure~\ref{fig:roadmap}), which we now numerically investigate in Section~\ref{sec:results}.

\section{Results}\label{sec:results}
\subsection{Comparison of the impact of rate-limited mechanisms}\label{sec:results_ratelim}

\begin{figure}
    \centering
\includegraphics[width=0.7\textwidth]{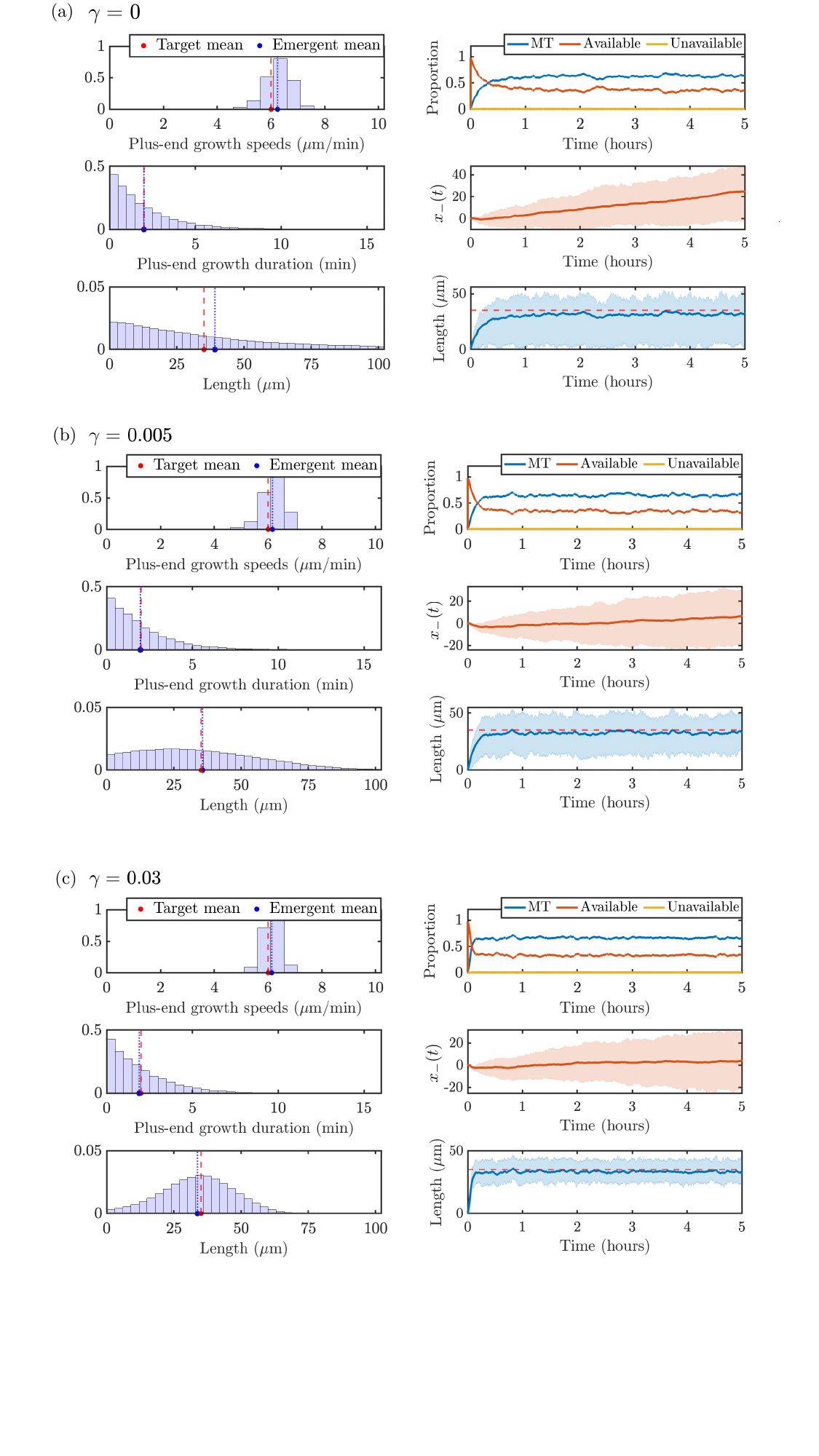}
    \caption{Target and emergent distributions of measurements and observables from the stochastic model for 10 trials with $N=20$ MTs, $T_{\mathrm{tot}} = 1000\mu$m and {\bf{(a)}} $\gamma = 0$, {\bf{(b)}} $\gamma = 0.005$, and {\bf{(c)}}, $\gamma = 0.03$. In each subpanel, {\bf{(Top left)}} plus-end growth speeds, {\bf{(middle left)}} plus-end growth durations, {\bf{(bottom left)}} microtubule lengths, {\bf{(top right)}} average tubulin pool allocation, {\bf{(middle right)}} minus-end position, and {\bf{(bottom right)}} average microtubule length. Line cloud represents the interquartile region. }\label{fig:dashboard_baseline_gammas}
\end{figure}

We investigate predictions of microtubule measurements based on stochastic simulation results in qualitatively different parameter regimes. As outlined in \S~\ref{sec:parameterization}, we wish to satisfy constraints (C1), (C2), and (C3) in order to have stochastic outputs match experimentally observable quantities (Figure~\ref{fig:roadmap}). We set the total amount of tubulin to $T_{\mathrm{tot}} = 1000 \mu$m and use the prescribed, experimentally-informed, and implied model parameters outlined in Table~\ref{tab:ode_params}. The mean target and experimentally observable quantities (plus-end growth speed, plus-end growth duration, MT length) for $T_{\mathrm{tot}} = 1000 \mu$m are shown with red dashed lines in the left panels of Figure~\ref{fig:dashboard_baseline_gammas} for the three values of the magnitudes of length-dependent catastrophe $\gamma$ investigated in Figure~\ref{fig:H-L-star}: (a) $\gamma = 0$, (b) $\gamma = 0.005$, and (c) $\gamma = 0.03$. We show plus-end growth duration ($1/\lambda_{g\to s,+}^*=2$~min) instead of the switching rate since it is a more interpretable quantity that can be obtained from experimental observations. The left panels of each parameter dashboard show that the distributions of the corresponding quantities that emerge from the stochastic model have averages that agree well with the target parameters (blue dotted lines). This provides evidence that the ODE model steady-state approximation provides useful insights for the stochastic model calibration. The right dashboard panels show time-course information on tubulin allocation, minus-end position, and microtubule length. Note that each curve corresponds to the mean quantity across 10 simulations, and the cloud denotes the interquartile spread of each output. 

The measurements used to generate the left panels of Figure~\ref{fig:dashboard_baseline_gammas} are extracted after $60$ minutes of simulation time pooled over the 10 simulations, to ensure that the system has achieved steady-state. The top left panels pool all strictly positive speeds of the MT plus ends from the stochastic simulation, to be consistent with experiments that track the velocities when the plus ends are actively moving. Turning on the length-dependent catastrophe mechanism shows a tightening of the emergent speed distribution (see Figure~\ref{fig:dashboard_baseline_gammas}b), due to the additional regulation of the microtubule lengths. The middle left panels show all the durations of growth at the plus ends of the microtubules, which are similar across the parameter ranges studied and show very good agreement with the target time spent in the growth state. The bottom left panels pool all strictly positive MT lengths. Interestingly, the length distribution for the $\gamma = 0$ case in Figure~\ref{fig:dashboard_baseline_gammas}a predicts an exponential-like MT length distribution, while adding on the length-dependent catastrophe mechanism in Figures~\ref{fig:dashboard_baseline_gammas}b and~\ref{fig:dashboard_baseline_gammas}c predicts a length distribution with a peak around the desired target length.

The right panels of Figure~\ref{fig:dashboard_baseline_gammas} illustrate emergent properties of the stochastic simulations through time, as introduced in Section~\ref{sec:model_obs}. The top right panels show the distribution of tubulin across available, unavailable, and microtubule pools. The length-dependent catastrophe mechanism limits the amount of tubulin in MTs and frees it up to the available pool, as expected. The middle right panels show the mean (red line) and the interquartile region (red region) of the MT minus end position trajectories, pooled across all MTs and all times. These panels show that, in the $\gamma=0$ regime (corresponding to no length-dependent catastrophe), we predict a larger overall movement of the entire microtubules in the dendrite segment. Finally, the bottom right panel shows the mean (blue solid line) and the interquartile region (blue cloud) of microtubule length over time. Compared to the length distribution in the bottom left panels in Figures \ref{fig:dashboard_baseline_gammas}a--c, the microtubule lengths over time include zero-length microtubules. Therefore, the average length as a function of time will be slightly smaller than the average length in the corresponding histogram. As $\gamma$ increases, the microtubule length distribution also tightens over time, where $\gamma = 0$ has the largest interquartile region.

\subsection{The predicted response to tubulin manipulation} \label{sec:results_tubmanip}

\begin{figure}
    \centering
    \includegraphics[width=0.8\textwidth]{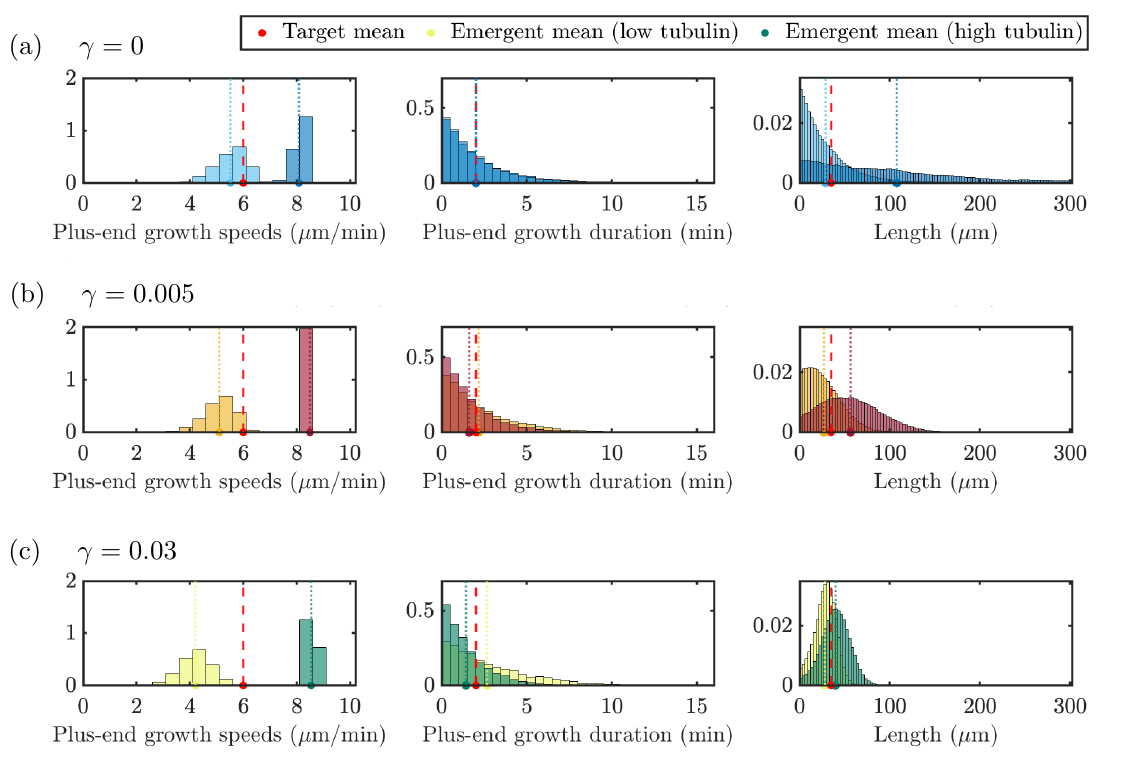}
    \caption{Distributions of target quantities for tubulin titration experiment with tubulin amounts $T_{\mathrm{tot}} = 700 \mu$m (light hues) and $T_{\mathrm{tot}} = 4000\mu$m (dark hues) for $N = 20$ MTs with input parameters set  for $\gamma = 0$ {\bf{(top, blue)}}, $\gamma=0.005$ {\bf{(middle, red)}} and $\gamma=0.03$ {\bf{(bottom, green)}}. Target and emergent mean values for various tubulin levels are denoted as in the dashboard figures.}
    \label{fig:tubulin_titration_hist}
\end{figure}

The previous section investigates how different length-regulating mechanisms impact the MT behavior in stochastic simulations with prescribed input parameters. The neuronal cell environment may however change through time and developmental stages, and in particular it may be characterized by different levels of tubulin. We are interested in predicting how these distinct tubulin levels may impact the MT growth and regulation. We assume that some growth characteristics of MTs would be fixed in varying tubulin environments. Namely, parameters $F_{1/2}, \bar{v}_{g}^{\pm}$, and $\bar{\tau}_g^{\pm}$ are likely innate characteristics of the cell that are unlikely to change as tubulin is varied. Therefore, we investigate a proposed tubulin titration experiment where we first parameterize our stochastic simulation at $T_{\mathrm{tot}} = 1000\mu$m based on the procedure in Section~\ref{sec:parameterization} and observe how emergent MT properties are affected by subsequently varying tubulin. For each level of tubulin, we perform 10 simulations and aggregate the emergent distributions and time-course data as described in the previous section. 

Figure~\ref{fig:tubulin_titration_hist} illustrates the impact of two extreme tubulin levels on the distributions of the emergent quantities and how the means of these distributions relate to the target parameters. Again, these distributions represent quantities taken after 60 minutes to ensure the system is at steady state. Microtubule length, plus-end growth speeds and plus-end growth durations distributions are shown for $T_{\mathrm{tot}} = 700\mu$m (light hue) and $T_{\mathrm{tot}} = 4000\mu$m (dark hue) with $N= 20$ microtubules. In the top panel, we show results for model simulations with no length-dependent catastrophe mechanism ($\gamma = 0$, blue). The emergent mean growth speed in the bottom left panel is below the target mean growth speed for the low tubulin level, and there is a wide spread of speeds observed in the simulation. We observe similar emergent distributions of the plus-end growth durations across tubulin levels, but also significantly different emergent length distributions for low and high tubulin. When tubulin is low, the emergent mean length is close to the target length, and most microtubules are shorter than the mean length. This is not the case for the high tubulin level, where the emergent mean length is much larger than the target length, and the exponential-like distribution shows much longer microtubules present in the system. 

For small length-dependent catastrophe ($\gamma = 0.005$, red) in the middle row of Figure~\ref{fig:tubulin_titration_hist}, the top right and center panels show that, for low tubulin, the target mean length and target plus-end growth duration are achieved, but the emergent growth velocity is below the target plus-end growth velocity. When tubulin is increased to $T_{\mathrm{tot}} = 4000\mu$m, the emergent growth velocity is larger than the target growth velocity and the distribution of speeds is much tighter. The emergent plus-end growth duration is also similar to the target duration, but we see slightly different emergent distributions of the plus-end growth duration. The emergent mean MT length in the top right panel is slightly larger than the target length and for both tubulin amounts, we find that the distribution of MT lengths is not exponential. When length-dependent catastrophe is increased further to $\gamma = 0.03$ (shown in green) in the bottom row of Figure~\ref{fig:tubulin_titration_hist}, there is little difference between the low and high tubulin distributions, and the distributions of MT lengths are symmetric about the target mean.

\begin{figure}
    \centering
      \includegraphics[width=\textwidth]{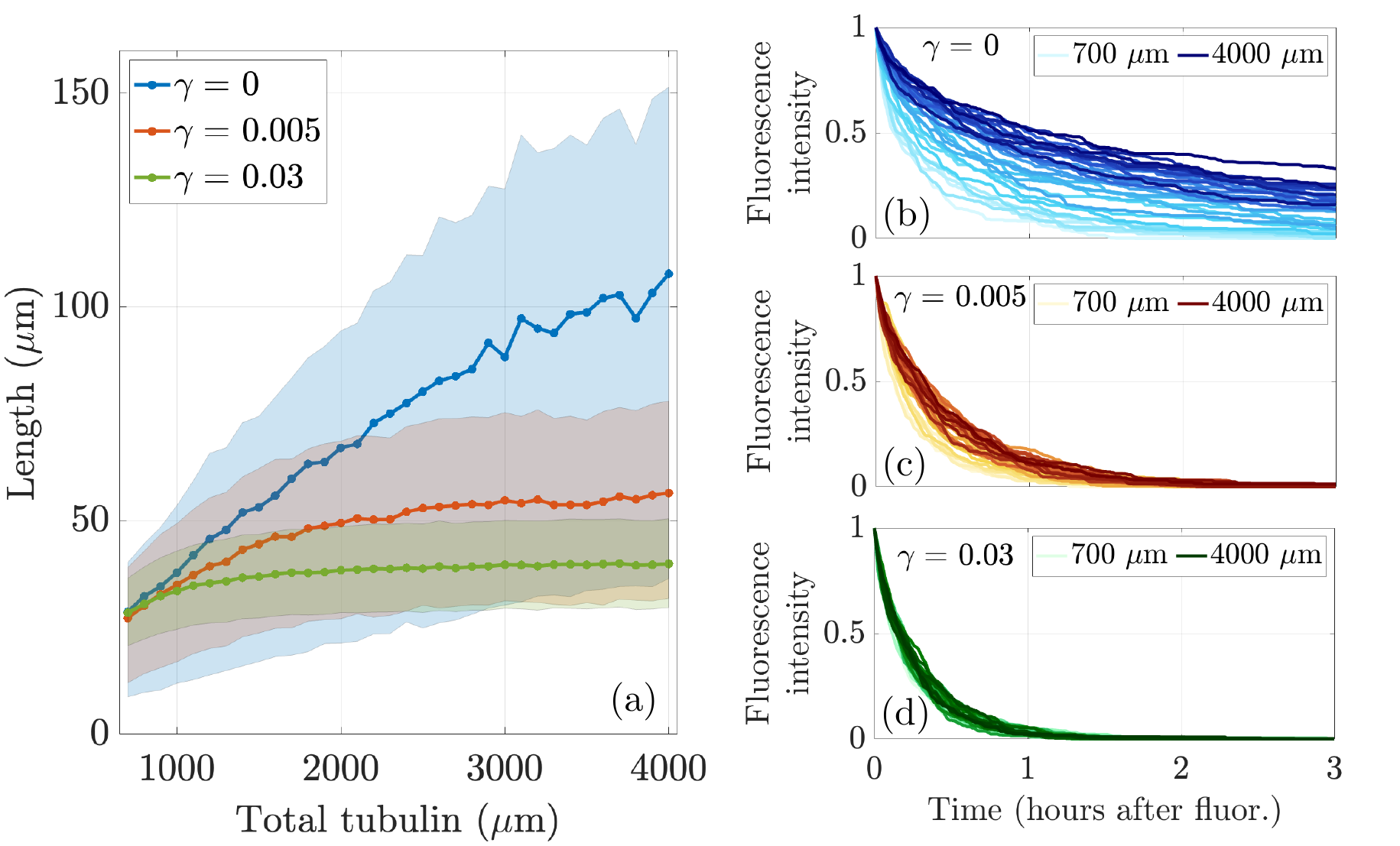}
    \caption{A contrast of responses to manipulating the level of tubulin in the cell for $N=20$ MTs when length-dependent switching to shrinking is present (red curves, $\gamma = 0.005$, green curves, $\gamma = 0.03$) and is not (blue curves, $\gamma = 0$). {\bf{(Left)}} Mean and interquartile region of average microtubule length as a function of total tubulin ($\mu$m).  {\bf{(Right)}} Average relative fluorescence intensity in models of tubulin photoconversion for various total tubulin levels. Colors in (b), (c), and (d) correspond to $\gamma = 0$, $\gamma = 0.005$, and $\gamma = 0.03$, respectively. }
    \label{fig:tubulin_titration}
\end{figure}

We now aim to understand how the calibration and prediction observables change as we vary tubulin across a range of levels, so that $T_{\mathrm{tot}} = \{700, 800, \dots, 3900, 4000\} \mu$m in the stochastic simulation with $N= 20$ microtubules and $L_{*} = 35 \mu$m. We show the emergent mean and interquartile range of microtubule length in the left panel of Figure \ref{fig:tubulin_titration} for 10 simulations each, where the tubulin-limited only case is shown in blue, the $\gamma = 0.005$ case is shown in red, and the $\gamma = 0.03$ case is shown in green. For low tubulin levels $(T_{\mathrm{tot}} \leq 1000 \mu$m), the $\gamma = 0$ and $\gamma = 0.005$ simulations result in similar mean microtubule lengths. As tubulin increases, the emergent mean MT length continues to increase for $\gamma = 0$, while for $\gamma=0.005$, the emergent mean MT length increases more slowly and appears to reach a plateau. The interquartile region for the $\gamma = 0$ case is much wider than the interquartile region for $\gamma=0.005$, indicating that there is a larger variance of microtubule lengths for the case with no length-dependent catastrophe. When $\gamma = 0.03$, length is very tightly regulated with a small interquartile region centered close to $35 \mu$m. Note that the mean MT lengths for $\gamma = 0, \gamma = 0.005,$ and $\gamma = 0.03$ found in the stochastic simulations exhibit similar behavior to the analytical results shown in Figure~\ref{fig:H-L-star} in Section~\ref{sec:ode}. 

Prediction observables also change as tubulin is varied. In the right panels of Figure \ref{fig:tubulin_titration}, we show the average fluorescence turnover curve at each tubulin amount for (b) $\gamma = 0$, (c) $\gamma = 0.005$, and (d) $\gamma = 0.03$. The darker colors correspond to larger tubulin levels for each $\gamma$. As tubulin increases, the average fluoresence turnover curve moves up and to the right for both choices of $\gamma$. This indicates that more fluorescent tubulin is persisting in microtubules for more time at large tubulin amounts compared to smaller tubulin amounts. Comparing Figure~\ref{fig:tubulin_titration}b, Figure~\ref{fig:tubulin_titration}c, and Figure~\ref{fig:tubulin_titration}d, we see a larger spread in fluorescence turnover for $\gamma = 0$ compared to $\gamma = 0.005$ and $\gamma = 0.03$. As tubulin increases for $\gamma = 0$, it takes longer for fluorescent tubulin to turnover, indicating that for large amounts of available tubulin, the protein stays in the microtubules at the end of the simulation. For some tubulin levels, the timescale of turnover is consistent with previous experiments \citep{Tao2016}. This behavior is not seen for $\gamma = 0.005$, where all tubulin levels exhibit a decay in fluorescence that results in no fluorescent tubulin left in the activation window at 3 hours. When $\gamma = 0.03$, the turnover is even more tightly regulated for all levels of tubulin. A direct comparison with actual experimental data is a subject of further investigation.

\section{Discussion}\label{sec:discussion}

Uncovering minimal mechanisms of MT length regulation in living cells is crucial for ultimately answering questions about cytoskeleton dynamics and organization in healthy as well as post-injury settings. For example, in healthy neurons, a stably organized microtubule cytoskeleton needs to persist for the lifetime of the animal. Similar microtubule lengths should be present for years or even decades, while individual microtubules grow, shrink and turn over. We therefore develop a model that can be used in concert with available data on MT turnover dynamics from \textit{in vivo} experiments. We consider a reduced differential equations model of MT dynamics at both ends using a traditional rate-balance compartmental model enhanced by terms that approximate intrinsically stochastic events like hitting the zero-length MT boundary. Steady-state analysis of this model informs our choices of parameters for a more extensive stochastic model of MT polymerization and depolymerization through switching between growth and shrinking states at both MT ends. The particular model parameterization we investigate here is motivated by experimental studies on MT turnover dynamics in segments of dendrites of \textit{Drosophila} neurons \citep{Feng2019, Rolls2021}. 

In contrast to prior studies, this modeling approach allows us to predict and match the spread of MT end speeds observed experimentally in the system inspiring this work. In addition, the model outputs predictions about the overall MT movement in the spatial domain, as well as about how tubulin protein may be allocated across MT and free tubulin pools \textit{in vivo}. It is worth noting that our parameterization approach requires determining the rates of switching from shrinking to growth at both MT ends based on the other available parameters. These rate estimates (and the associated run lengths spent shrinking) from both the deterministic and the stochastic models are different from those estimated experimentally. This could be due to the fact that run lengths may be cut off in experiments, in which data was captured for at most 20~minutes. This could either occur because of the MT switching back to a growth phase, or because the MT has reached zero length.

The most prominent contrast we draw in our modeling framework is between two {physically plausible} mechanisms of MT length control: tubulin-limitation and length-dependent catastrophe. We find that both mechanisms can yield MT length distributions with any given mean, but length-dependent catastrophe gives rise to much more variation in distribution shape. In the tubulin-limited regime, MT lengths appear essentially exponential in distribution, while in the length-dependent catastrophe regime, we observe anywhere from monotonically decreasing exponential-like distributions to unimodal Gaussian-looking distributions. Interestingly, this is consistent with experimental observations that have suggested that the steady-state length distribution of MTs is either exponential or Gaussian \citep{cassimeris1986dynamics,schulze1986microtubule,verde1990regulation,gliksman1992okadaic,gardner2011depolymerizing,gardner2013microtubule}. {While the limited tubulin availability in neurons is not controversial, the exact molecular basis for the length-dependent catastrophe mechanism is less well established. Based on the modeling results presented here (in Sections~\ref{sec:results_ratelim},\ref{sec:results_tubmanip}), Gaussian-like MT length distributions with narrow peaks seem to suggest that both length-dependent catastrophe and tubulin limitations contribute to MT length control. Wider distributions would imply that tubulin availability may not be as large of a limitation to MT growth. The model provides a useful tool for helping distinguish between contributions of these mechanisms, and invites further study when experimental data on MT lengths \textit{in vivo} becomes available.}

Our stochastic model also predicts measurements that could be potentially observable in the future, such as the effective direction of movement of the MT ends. While the emergent dynamics of the minus ends requires further study, such insights on overall MT movement are typically not available from prior models. These insights can be further studied in the context of recent \textit{in vivo} studies such as \citep{Feng2019}, which found that minus ends display sustained growth and that these dynamics are important for getting minus-end-out MTs into distal regions of healthy developed dendrites, as well as into developing and regenerating dendrites. 

Similarly, the tubulin titration simulations provide a useful thought experiment for the impact of tubulin availability on MT regulation in the changing cell environment. In Section~\ref{sec:results_tubmanip}, we illustrate how low and high tubulin amounts influence the predicted growth speeds of MT ends as well as the MT lengths themselves. The two mechanisms of length regulation considered also generate distinguishable distributions of MT lengths and significant changes in the decay timescales for the model fluorescent turnover experiments. 

The model thus generates testable predictions for experiments that manipulate tubulin in neurons, which have been proposed and could be used to validate hypothesized model mechanisms. {In particular, if future experiments allow for control of tubulin levels in living cells, our model predicts that stabilization of mean and variance of the MT lengths with increasing tubulin levels would suggest a significant contribution of the length-dependent catastrophe mechanism to MT length regulation. If photoconversion experiments are feasible at varying tubulin amounts, our model suggests that abundant tubulin regimes would be characterized by large changes in the timescale of fluorescent tubulin turnover with increasing tubulin levels (see results in Section~\ref{sec:results_tubmanip}).}

A limitation of our current stochastic model is that it considers a fixed number of MTs, so that complete shrinking of a MT is followed by nucleation, or seeding, of a new MT. This modeling choice is motivated by the fact that nucleation is known to be heavily controlled in the cell. In particular, nucleation is heavily regulated in neurons \citep{chen2012axon,hertzler2020kinetochore}, and plays a role in axon injury responses \citep{chen2012axon} and dendrite regeneration \citep{nye2020receptor}. However, it is not clear how robust this regulation is with respect to major changes in tubulin availability or other potential experimental manipulations. Depending on forthcoming experimental observations, allowing for a variable number of MTs could be an important future direction for enhancing our model. The proposed stochastic model and the parameterization procedure based on steady-state analysis of the reduced deterministic model can also extend to \textit{in vivo} experiments in other animal systems, which are likely to be characterized by different MT lengths and polymerization speeds.

\bmhead{Acknowledgments}
This work was supported by NIH grant R01NS121245. ACN was partially supported by NSF grant DMS-2038056.

\bmhead{Data availability}
Experimental data is referenced appropriately throughout the manuscript and all simulation data is available upon request.

\begin{appendices}
\section{Stochastic model implementation details}
\label{app:stochmodel}

\begin{figure}[h]
    \centering
    \includegraphics[width=0.7\textwidth]{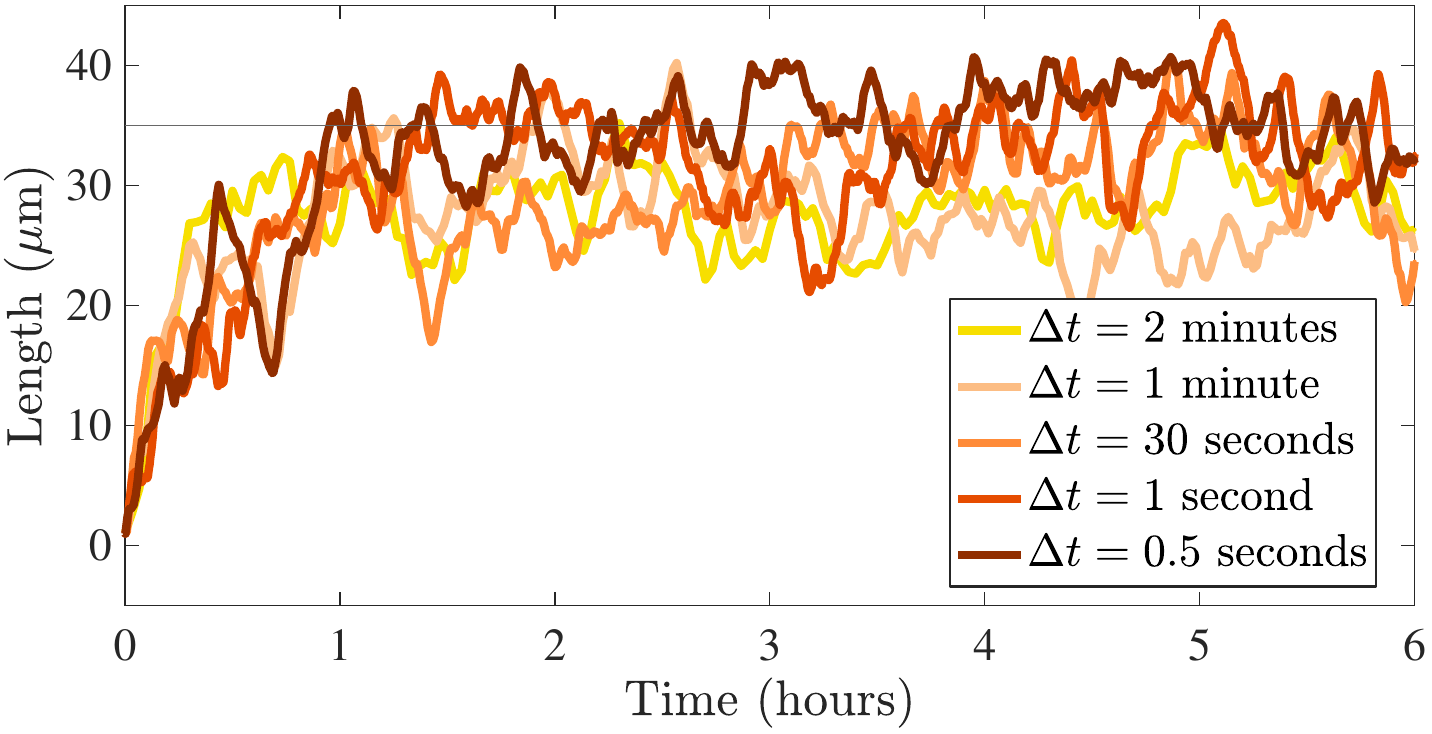}
    \caption{Average MT length over time for different stochastic time steps, $\Delta t$ for $\gamma = 0$, $L_* = 35\mu$m, and $N=20$ MTs for one realization.  }
    \label{fig:time_titration}
\end{figure}

To implement the stochastic model of MT polymerization and depolymerization in Section~\ref{sec:stochastic}, we {implement 
an Euler-based stochastic scheme that is time- rather than event-driven, and thus we update events at fixed time steps $\Delta t$. Similar to the leap condition in Gillespie's tau-leaping scheme \citep{gillespie2001approximate,gillespie2007stochastic}, we choose a time interval $\Delta t$ that is small enough so that reaction propensities are relatively fixed in time interval $[t_i, t_i+\Delta t = t_{i+1}]$, and large enough so that many reactions are still happening in each step. Since we are interested in the emergent behavior of dynamic MT instability, our approach differs from a strict implementation of tau-leaping in that we avoid resolving tubulin unit scale reactions and instead focus on MT end state switching reactions. Tau-leaping methods face significant challenges when one or more of the state variables are near boundary values. In our case, the pool of available tubulin is often near zero. This also motivates our choice of time scale, and our decision to treat tubulin resources as a continuum, rather that in discrete units.} To determine how small time step $\Delta t$ should be, we performed an experiment shown in Figure~\ref{fig:time_titration} to see the effect of time step size on average length. We see MT behavior is qualitatively the same for time step sizes less than or equal to one second, so we utilize $\Delta t = 1 \text{ second}$.

To determine how long each end spends in a certain state during time step $\Delta t$, we first update the length-dependent catastrophe rate based on the length of the microtubules at the previous time step, $t_i$, denoted as $L_i = L(t_i)$. Therefore, the catastrophe rate at $t_{i+1} = t_i + \Delta t$ will take the following form 
\begin{equation}\label{eq:app_lengthdepcat}
   \lambda^{\pm}_{g\rightarrow s}(L_i) =   \max\left(\lambda_{\min}, \lambdags^\pm + \gamma L_0 \phi\left(\frac{L_i}{L_0}\right)\right),
\end{equation}
where $\catform(x) = x - 1$. Again, for small $L$ and large $\gamma$, it is possible that $ \lambda_{min} > \lambdags^\pm + \gamma L_0(\frac{L(t)}{L_0} - 1)$, so then the catastrophe rate will be $\lambda_{min}$ for short MTs. The updated catastrophe rates are then used to draw times spent in growth phase.

\begin{figure}
    \centering
    \includegraphics[width=0.4\textwidth]{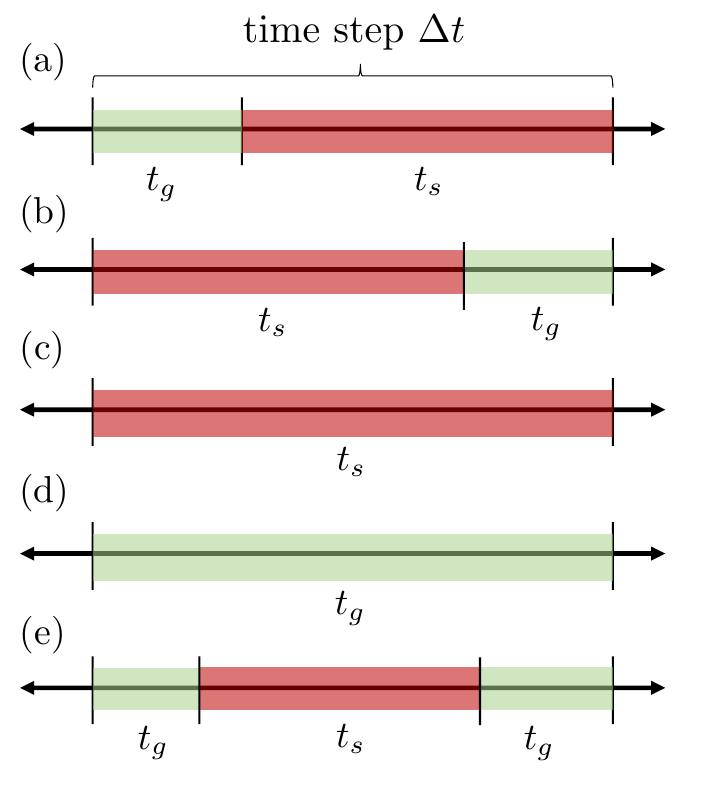}
    \caption{Sketch of admissible growth and shrinking times ($t_g$ and $t_s$, respectively) within a single time step $\Delta t$ of the stochastic simulation. Green color denotes growth and red denotes shrinking. In a single $\Delta t$, we only permit one shrinking phase to occur. }
    \label{fig:times}
\end{figure}

At each time step from $t_i$ to $t_i + \Delta t$, each microtubule end is in either a shrinking or growth state. For each end of the microtubule, the time spent in growth and shrinking, denoted as $t_g$ and $t_s$, respectively, are given by random variables drawn from the appropriate exponential distribution, where
\begin{equation}
    t_g^{i+1} \sim \text{Exp}(\lambda_{g\rightarrow s}(L_i)), \qquad  t_s^{i+1} \sim \text{Exp}(\lambda_{s\rightarrow g}).
    \label{eq:drawtimes}
\end{equation} 
We then determine if the microtubule exits the time step in either growth or shrinking based on the idea that only one shrinking phase can occur in a time step $\Delta t$ (see Figure~\ref{fig:times}). For example, if a microtubule end entered in growth and $t_g > \Delta t$, then the microtubule would also exit the time step in growth (see Figure \ref{fig:times}d). Similarily, if a microtubule end entered in growth and $t_g + t_s < \Delta t$, then the microtubule also exits in growth, shown in Figure \ref{fig:times}e. 

Once the time spent in each state has been determined, we then perform the shrinking events, which involve updating the microtubule lengths and adding the lost tubulin to unavailable tubulin bin. If the shrinking events resulted in a completely-catastrophed microtubule (that is, a microtubule with length less than or equal to zero), we reseed the microtubule at its previous position with length zero and both ends in a growth state. This reseeding procedure ensures that the number of microtubules remains the same throughout the simulation. 

After performing shrinking events, we update the available tubulin pool from the unavailable tubulin pool at rate $\tau_{tub}$ and then implement growth events, which depend on the available tubulin pool. To allow for a smooth transition in growtexth speed as the available tubulin pool becomes very small, we model the dependence of the microtubule growth on the available pool using Michaelis--Menten kinetics. For example, suppose the desired amount of tubulin for growth for a single MT is $v_{g} t_g$, where $v_{g}$ is the growth velocity and $t_g$ is the time spent in growth for that MT in $\Delta t$. Since growth is dependent on tubulin availability, we update the desired growth length for the $j$-th MT to be
\begin{equation}
    \tilde{x}_j = \frac{F}{F_{1/2} + F}v_g t_g, 
\end{equation}
where $F$ is the available tubulin amount and $F_{1/2}$ is the same Michaelis--Menten constant from Equation~\eqref{eq:F-half}.

We implement the growth events for all MTs in $\Delta t$ and decrease $F$ based on the total desired polymerization length,  $\sum_{i=1}^N \tilde{x}_i$. If tubulin resources are limited, it is possible that the available tubulin amount is less than the total desired polymerization length, so that $F < \sum_{i=1}^N \tilde{x}_i$. To account for this, we divide the available tubulin  amount proportionally to each MT based on the desired polymerization length for that time step. Therefore, the final length grown for the $j$-th MT, $x_j$, is  
\begin{equation}
    x_j = \frac{\tilde{x}_j}{\sum_{i=1}^N \tilde{x}_i} F \,. 
\end{equation}
We then deplete the tubulin used to grow MTs from the available tubulin pool, $F$, so that 
\begin{equation}
    F(t_{i+1}) = F(t_i) - \sum_{i = 1}^N x_i \,.
\end{equation}

\end{appendices}

\bibliography{MTRef}

\end{document}